\documentclass[preprint,11pt]{elsarticle}

\makeatletter
\def\ps@pprintTitle{%
 \let\@oddhead\@empty
 \let\@evenhead\@empty
 \def\@oddfoot{\centerline{\thepage}}%
 \let\@evenfoot\@oddfoot}
\makeatother

\usepackage{amssymb}
\usepackage{amsthm}
\usepackage{amsmath}
\usepackage{algorithm}
\usepackage{algpseudocode}
\usepackage{graphicx}
\usepackage{subcaption}
\usepackage{adjustbox}
\usepackage{empheq}
\usepackage{ulem}

\usepackage{tikz,xcolor}
\usepackage{pgfplots}
\pgfplotsset{compat=newest}
\usepackage{mathtools}
\usepackage{mathrsfs}

\newtheorem{theorem}{Theorem}
\newtheorem{remark}{Remark}

\usepackage{stmaryrd}
\usepackage{tikz}
\usetikzlibrary{positioning}
\usetikzlibrary{shapes,arrows}

\biboptions{sort&compress}

\usepackage{listings}
\usepackage{color} 
\definecolor{mygreen}{RGB}{28,172,0} 
\definecolor{mylilas}{RGB}{170,55,241}

\newcommand{\R}{\mathbb{R}}
\newcommand{\eps}{\varepsilon}

\begin{document}

\begin{frontmatter}



\title{Variational Formulations for Explicit Runge-Kutta Methods}


\author[1]{Judit Mu\~noz-Matute}
\author[1,2,3]{David Pardo}
\author[4,5]{Victor M. Calo}
\author[1]{\\Elisabete Alberdi}

\address[1]{University of the Basque Country (UPV/EHU), Leioa, Spain}
\address[2]{BCAM-Basque Center for Applied Mathematics, Bilbao, Spain}
\address[3]{IKERBASQUE, Basque Foundation for Science, Bilbao, Spain}
\address[4]{Applied Geology, Western Australian School of Mines, Faculty of Science and Engineering, Curtin University, Perth, WA, Australia 6845}
\address[5]{Mineral Resources, Commonwealth Scientific and Industrial Research Organisation (CSIRO), Kensington, WA, Australia 6152}

\begin{abstract}
Variational space-time formulations for Partial Differential Equations have been of great interest in the last decades. While it is known that implicit time marching schemes have variational structure, the Galerkin formulation of explicit methods in time remains elusive. In this work, we prove that the explicit Runge-Kutta methods can be expressed as discontinuous Petrov-Galerkin methods both in space and time. We build trial and test spaces for the linear diffusion equation that lead to one, two, and general stage explicit Runge-Kutta methods. This approach enables us to design explicit time-domain (goal-oriented) adaptive algorithms. 
\end{abstract}

\end{frontmatter}


\section{Introduction}\label{S1}
Adaptive algorithms for Partial Differential Equations (PDEs) are powerful tools to produce optimal grids that seek to minimize the computational cost. For time dependent problems, it is common to employ time-marching schemes and adapt the time step size and/or the spatial mesh size employing a posteriori error estimates \cite{chen2004adaptive,nicaise2005posteriori,collier2011diffusive}. On the other hand, there exist adaptive strategies based on space-time Finite Element Methods (FEM). In some algorithms, the space-time tensor-product of the trial and test functions is assumed. Then, employing discontinuous-in-time basis functions, the space-time FEM can be reinterpreted as a time-marching scheme \cite{schotzau2000hp,schieweck2010stable}. Here, the approximation orders (in space and time) as well as the mesh size and the time step size can be adapted \cite{schotzau2000time,werder2001hp}. But the main benefit of the space-time FEM is that they can be employed to build unstructured space-time meshes \cite{erickson2005building,abedi2006h,abedi2006space,miller2008spacetime}. 

To perform goal-oriented adaptivity \cite{bangerth2001adaptive,schmich2008adaptivity,besier2012goal,csimcsek2015duality}, we require a space-time variational formulation of the problem. In this kind of algorithms, we represent the error in the quantity of interest as an integral over the whole space-time domain that is subsequently expressed as a sum of local element contributions, which we use for adaptivity. A full space-time variational formulation allows such representation \cite{diez2007goal}. However, all existing space-time variational formulations for PDEs lead to implicit methods in time when they are thought as time marching schemes. For that reason, adaptive strategies based in space-time FEM as well as time-domain goal-oriented adaptive processes are performed employing implicit methods in time. Our focus is to design goal-oriented adaptive algorithms employing explicit methods in time since in many instances they can be computationally cheaper than implicit ones.

In the resolution of time-dependent PDEs, we commonly discretize independently the spatial and temporal variables (also called semidiscretization or Method of Lines) \cite{schiesser2012numerical}. First, the spatial variable is discretized by the Finite Element Method (FEM) to obtain a system of ODEs. The resulting system is subsequently solved employing time stepping schemes. The alternative idea of using variational space-time methods was well established in the late eighties and early nineties \cite{eriksson1987error,aziz1989continuous,french1993space,johnson1993discontinuous}. Hughes and Hulbert \cite{hughes1988space,hulbert1990space} proposed a stabilized space-time FEM for hyperbolic problems. They showed that the oscillations present in the solutions were considerably reduced by employing space-time variational formulations rather than using semidiscretizations. Nowadays, it is well known that some low-order space-time FEM are algebraically equivalent to some semidiscretizations \cite{bangerth2010adaptive}. For example, the discontinuous Galerkin method using constants in time (usually denoted by dG(0)) leads to the Backward Euler method. The continuous Petrov-Galerkin method in time with linear trial and constant test functions (denoted cGP(1)) is equivalent to the Crank-Nicholson method. Recently, higher order dG(k)- and cGP(k)-methods have been developed and analyzed for parabolic and hyperbolic problems \cite{hussain2011higher,kocher2014variational,ahmed2015adaptive,ern2017guaranteed}.

There also exist variational formulations of time marching schemes in the context of Ordinary Differential Equations (ODEs) \cite{bottasso1997new}. Delfour et. al. \cite{delfour1981discontinuous,delfour1986discontinuous} and Hulme \cite{hulme1972one} showed that it is possible to obtain classical schemes like Runge-Kutta methods by employing Galerkin methods for initial value problems together with quadrature formulas. Estep and French \cite{estep1994global,estep1995posteriori} derived error bounds for the continuous and discontinuous Galerkin methods to efficiently adapt the time step size. More recently, Estep and Stuart studied the dynamical behavior of discontinuous Galerkin methods for ODEs in \cite{estep2002dynamical}; in \cite{collins2015posteriori}, Collins et. al. derived an a posteriori error estimation of explicit schemes; and Tang et. al. provided in \cite{tang2012time} a unified framework of finite element methods in time.

In this work, we construct a variational formulation of explicit Runge-Kutta methods for parabolic problems.
Such formulation could be applied to design explicit-in-time (and consequently cheaper) goal-oriented adaptive algorithms, to build new time-stepping schemes and also to extend the existing space discretizations like IGA \cite{hughes2005isogeometric,gomez2008isogeometric,vignal2017energy}, DPG \cite{demkowicz2014overview} and Trefftz \cite{egger2015space}, to time domain problems.

First, we derive the Forward Euler method (one stage explicit Runge-Kutta). For trial functions, we build a family of piecewise linear polynomials, depending on a parameter $\eps$, that are globally-continuous in time. We construct these functions in such a way that when we take the limit $\eps\to0$, they become discontinuous and piecewise constant in time. We prove that by employing these trial functions and selecting piecewise constant test functions, we recover the Forward Euler method. This construction supports dynamic meshes in space, i.e, we allow different spatial discretizations per time interval. In order to obtain square mass matrices, we displace in time the spatial discrete spaces of the test space with respect to the trial space. This displacement leads to a Petrov-Galerkin method both in space and time. 

For a general number of Runge-Kutta stages, we perform a similar construction. That is, we generate the trial space by a family of piecewise polynomials that are globally continuous in time and depend on a parameter $\eps$. For the test space, we define a family of piecewise polynomials that are discontinuous in time. By substituting them into the variational formulation and treating the coefficients of the polynomials as unknowns, we obtain integrals of the products of trial and test functions. Then, we establish some conditions that these time integrals must satisfy. First, we state the necessary orthogonality conditions needed to obtain an explicit method. We also define non-orthogonality conditions by matching the remaining integrals with the entries of the Butcher's table that define the Runge-Kutta methods. Finally, taking the limit $\eps\to0$ and performing analytic integration, we obtain a system of nonlinear equations. By solving this system, we obtain the coefficients of the trial and test functions for any stage Runge-Kutta method. However, for a large number of stages (above 5) the system becomes hard to solve. 

This article is organized as follows. Section \ref{S2} describes the strong and weak formulations of the linear diffusion equation we use to develop the theory. In Section \ref{S3}, we derive a discontinuous-in-time Petrov-Galerkin formulation of the one, two, and general $s$-stage Runge-Kutta methods, providing some examples. Section \ref{S4} explains the conclusions and the possible extensions of this work. Finally, in \ref{A1}, we express in matrix form the nonlinear system of equations we need to solve to obtain any explicit Runge-Kutta method. \ref{A2} provides a MATLAB code to solve it.

\section{Model Problem}\label{S2}
In this section, we state both the strong and weak formulations of the model problem we employ to develop the discontinuous-in-time Petrov-Galerkin formulations.
\subsection{Strong formulation}
Let $\Omega\subset\R^{d}$, where $d\in\{1,2,3\}$, and $I=(0,T]\subset\R$. We consider the \textit{linear heat (diffusion) equation}

\begin{equation}\label{heat}
\displaystyle{ \left\{
\begin{split}
u_{t}-\Delta u=&\;f&\mbox{in}&\;\Omega\times I,\\
u=&\;0&\mbox{on}&\;\partial\Omega\times I,\\
u(0)=&\;u_{0}\;&\mbox{in}&\;\Omega,\\
\end{split}
\right.} 
\end{equation}
where $u_{t}:=\partial u/\partial t$, $\Delta u = div(grad(u))$ is the Laplacian of $u$ and $\partial\Omega$ denotes the boundary of the spatial domain $\Omega$. 
The solution $u(\mathbf{x},t)$ of (\ref{heat}) could represent the temperature distribution in a body.
The source term $f(\mathbf{x},t)$ and the initial temperature distribution $u_{0}(\mathbf{x})$ are given data. For arbitrary Dirichlet (geometric) boundary conditions, we can modify the source term accordingly, thus making (\ref{heat}) a general statement.

\subsection{Weak formulation}
In order to obtain the weak formulation of (\ref{heat}), we multiply the diffusion equation by test functions $v$ of a suitable space $\mathcal{V}$ and we integrate over the whole domain $\Omega\times I$ 
$$\int_{I}\int_{\Omega}\bigg(vu_{t}-v\Delta u\bigg)\;d\Omega\;dt=\int_{I}\int_{\Omega}vf\;d\Omega\;dt, \;\forall v\in \mathcal{V}.$$
Integrating by parts in space the diffusion term and selecting test functions vanishing on $\partial\Omega$ we obtain
\begin{equation}\label{integ}
\int_{I}\int_{\Omega}\bigg(vu_{t}+\nabla v\cdot\nabla u\bigg)\;d\Omega\;dt=\int_{I}\int_{\Omega}vf\;d\Omega\;dt, \;\forall v\in \mathcal{V}.
\end{equation}
A sufficient condition for the above integrals to make sense is if all factors in the above products are in $L^{2}$ in both space and time.
For the space integrals, taking into account the diffusion term in (\ref{integ}), it seems natural that $u$ and $v$ should be in
$$V:=H_{0}^{1}(\Omega):=\{u\in L^{2}(\Omega)\;|\;\nabla u\in L^{2}(\Omega),\;u=0\;\mbox{on}\;\partial\Omega\},$$
and therefore, to guarantee integrability of the weak formulation, $u_{t}$ and $f$ should belong to $V':=H^{-1}(\Omega)$, which is the dual space of $V.$ In time, it is enough to ensure that all the functions in (\ref{integ}) are in $L^{2}.$

In the remaining of this article, we omit the spatial dependence of the functions, i.e., we write $u(t)$ instead of $u(\mathbf{x},t)$. We consider $u(t)$ and $v(t)$ as functions in time that take values in suitable Hilbert spaces \cite{salsa2008partial}, so in view of the sufficient integrability conditions discussed in the previous paragraph, we construct the following test space
$$\mathcal{V}:=L^{2}(I;V)=\left\{u:I\longrightarrow V\;|\;u\;\mbox{is}\;V\mbox{-measurable and}\;\int_{I}||u(t)||^{2}_{_{V}}\;dt<+\infty\right\},$$
which is the space of all integrable functions in time that take values in $V$. 

On the other hand, for the solution, we need $u\in\mathcal{V}$ and $u_{t}\in\mathcal{V'}:=L^{2}(I;V')$. Treating the initial solution as a Dirichlet condition in time, we define the following trial space
$$\mathcal{U}:=\{u\in \mathcal{V}\;|\;u_{t}\in \mathcal{V'},\;u(0)=u_{0}\}.$$ 
Finally, assuming that $f\in\mathcal{V'}$ and $u_{0}\in V$, the weak solution of problem (\ref{heat}) satisfies
\begin{equation}\label{varspacetime}
\int_{I}\left<v,u_{t}\right>dt+\int_{I}(\nabla v,\nabla u)dt=\int_{I}\left<v,f\right>dt,\;\forall v\in\mathcal{V},\\
\end{equation}
where $\left<\cdot,\cdot\right>$ denotes the duality pairing between the spaces $V$ and $V'$, and $(\cdot,\cdot)$ denotes the inner product in $L^{2}(\Omega)$. Now, denoting by
$$B^{t}(v,u):=\int_{I}\left<v,u_{t}\right>dt,\;\;\;B^{g}(v,u):=\int_{I}(\nabla v,\nabla u)dt,$$
\begin{equation}\nonumber
\begin{split}
B(v,u)&:=B^{t}(v,u)+B^{g}(v,u),\\
F(v)&:=\int_{I}\left<v,f\right>dt,
\end{split}
\end{equation}
we obtain the \textit{weak formulation} of problem (\ref{heat})
\begin{equation}\label{primal}
\displaystyle{ \left\|
\begin{tabular}{lll}
Find $u\in\mathcal{U}$ such that\\
$B(v,u)=F(v),\;\;\forall v\in\mathcal{V}$,\\
\end{tabular}
\right.} 
\end{equation}
where $F(\cdot)$ is a linear form and $B(\cdot,\cdot)$ is a bilinear form.

\subsection{Abstract discretization}
Selecting discrete space-time trial and test subspaces $\mathcal{U}_{\tau h}^{\eps}\subset\mathcal{U}$, $\mathcal{V}_{\tau h}\subset\mathcal{V}$, we define the fully discrete problem
\begin{equation}\label{primaldis}
\displaystyle{ \left\|
\begin{tabular}{lll}
Find $u_{\tau h}\in\mathcal{U}_{\tau h}^{\eps}$ such that\\
$B(v_{\tau h},u_{\tau h})=F(v_{\tau h}),\;\;\forall v_{\tau h}\in\mathcal{V}_{\tau h}$.\\
\end{tabular}
\right.} 
\end{equation}
In the following sections, we show that by selecting different trial and test subspaces in (\ref{primaldis}), we obtain different methods in time.

\section{Explicit Runge-Kutta methods}\label{S3}
In this section, we derive \textit{discontinuous Petrov-Galerkin} formulations that are equivalent to explicit Runge-Kutta methods in time. In space, we employ either a \textit{Spectral Element Method (SEM)} or a \textit{Finite Element Method (FEM)}.

We perform a partition of the time interval $\bar{I}=[0,T]$ as
$$0=t_{0}<t_{1}<\ldots<t_{m-1}<t_{m}=T,$$
and denote $I_{k}:=(t_{k-1},t_{k})$, $\tau_{k}:=t_{k}-t_{k-1},\;\forall k=1,\ldots,m$ and $\tau:=\max\limits_{1 \leq k \leq m}\tau_{k}$.

We allow dynamic meshes in space so we define a finite dimensional subspace of $V$ for each time step, i.e., $V_{h}^{k}\subset V,\;\forall k=0,\ldots,m,$ where $h$ is the largest element diameter of each dynamic mesh.

\subsection{Forward Euler method}\label{S3.1}
\begin{figure}[h]
\centering
\parbox{4cm}{
\begin{tikzpicture}
	\draw[->]	(0,0) -- (2.5,0) node[anchor=south] {\scriptsize{$t$}};
	\draw[->]	(0,0) -- (0,2.5) node[anchor=east] {};
	\draw	(-0.1,2) node{{\scriptsize 1}};
	\draw	(0,0) node[anchor=north] {\scriptsize{$t_{0}$}};
	\draw	(2,0) node[anchor=north] {\scriptsize{$t_{1}$}};
	\draw (2,2pt) -- (2,-2pt) node[anchor=north] {};
	\draw[thick] (0,2) -- (1.5,2);
	\draw[thick] (1.5,2) -- (2,0);
	\draw[thick] (2,0) -- (2.5,0);
	\draw[dotted] (2,0) -- (2,2);
	\draw[dotted] (1.5,2) -- (2,2);
	\draw[thick,red] (1.5,0) -- (2,0);
	\draw (1.75,0.2) node{{\textcolor{red}{\scriptsize{$\eps$}}}};
	\draw (0.9,2.4) node {\small{$\phi_{\eps}^{0}$}};
\end{tikzpicture}
}\parbox{5.4cm}{
\begin{tikzpicture}
	\draw[->] (0,0) -- (4,0) node[anchor=south] {\scriptsize{$t$}};
	\draw[->] (0,0) -- (0,2.5) node[anchor=east] {};
	\draw (-0.1,2) node{{\scriptsize 1}};
	\draw (0.5,0) node[anchor=north] {\scriptsize{$t_{k-2}$}};
	\draw (2,0) node[anchor=north] {\scriptsize{$t_{k-1}$}};
	\draw (3.5,0) node[anchor=north] {\scriptsize{$t_{k}$}};
	\draw (3.5,2pt) -- (3.5,-2pt) node[anchor=north] {};
	\draw (0.5,2pt) -- (0.5,-2pt) node[anchor=north] {};
	\draw (2,2pt) -- (2,-2pt) node[anchor=north] {};
	\draw[thick] (1.5,0) -- (2,2) --(3,2)-- (3.5,0);
	\draw[thick] (0,0) -- (1.5,0);
	\draw[thick] (3.5,0) -- (4,0);
	\draw[thick,red] (1.5,0) -- (2,0);
	\draw (1.75,0.2) node{{\textcolor{red}{\scriptsize{$\eps$}}}};
	\draw[thick,red] (3,0) -- (3.5,0);
	\draw (3.25,0.2) node{{\textcolor{red}{\scriptsize{$\eps$}}}};
	\draw[dotted] (2,0) -- (2,2);
	\draw[dotted] (0,2) -- (2,2);
	\draw[dotted] (3,2) -- (3.5,2);
	\draw[dotted] (3.5,0) -- (3.5,2);
	\draw (2.8,2.4) node {\small{$\phi_{\eps}^{k-1}$}};
\end{tikzpicture}
}\parbox{2.2cm}{ 
\begin{tikzpicture}
	\draw[->] (0,0) -- (2.5,0) node[anchor=south] {\scriptsize{$t$}};
	\draw[->] (0,0) -- (0,2.5) node[anchor=east] {};
	\draw (-0.1,2) node{{\scriptsize 1}};
	\draw (0.5,0) node[anchor=north] {\scriptsize{$t_{m-1}$}};
	\draw (2,0) node[anchor=north] {\scriptsize{$t_{m}$}};
	\draw (2,2pt) -- (2,-2pt) node[anchor=north] {};
	\draw (0.5,2pt) -- (0.5,-2pt) node[anchor=north] {};
	\draw[thick] (0,0) -- (1.5,0);
	\draw[thick] (1.5,0) -- (2,2);
	\draw[thick,red] (1.5,0) -- (2,0);
	\draw (1.75,0.2) node{{\textcolor{red}{\scriptsize{$\eps$}}}};
	\draw[dotted] (2,0) -- (2,2);
	\draw[dotted] (0,2) -- (2,2);
	\draw (2.3,2.3) node {\small{$\phi_{\eps}^{m}$}}; 
\end{tikzpicture}
}
\caption{Trial functions for the forward Euler method.}
\label{TrialFE}
\end{figure}
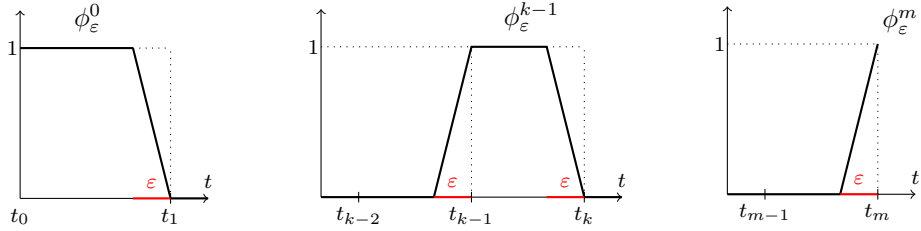

To appropriately define the discrete trial and test subspaces of (\ref{primaldis}) for this particular problem, we proceed as follows
\begin{equation}\label{trialtestFE}
\begin{split}
\mathcal{U}_{\tau h}^{\eps}&:=u_{0}\phi_{\eps}^{0}(t)+\mbox{span}\{u_{h}^{k}\phi_{\eps}^{k}(t),\;u_{h}^{k}\in V_{h}^{k},\;\forall k=1,\ldots,m\}\subset\mathcal{U},\\
\mathcal{V}_{\tau h}&:=\mbox{span}\{v_{h}^{k}\varphi^{k}(t),\;v_{h}^{k}\in V_{h}^{k},\;\forall k=1,\ldots,m\}\subset\mathcal{V},
\end{split}
\end{equation}
where the trial functions for $0<\eps<\min\limits_{1 \leq k \leq m} \tau_{k}$ are defined as (see Figure \ref{TrialFE})
$$\phi_{\eps}^{k-1}(t)=\begin{cases}
              \displaystyle{\frac{t-(t_{k-1}-\eps)}{\eps}},&t\in(t_{k-1}-\eps,t_{k-1}),\\ 
              1,&t\in[t_{k-1},t_{k}-\eps],\\
              \displaystyle{\frac{t_{k}-t}{\eps}},&t\in(t_{k}-\eps,t_{k}),\\
              0,&\mbox{elsewhere},\\
	      \end{cases}\;\;\;\forall k=2,\ldots,m,$$
\begin{minipage}[t]{6cm}
$\phi_{\eps}^{0}(t)=\begin{cases}
	        1,&t\in(t_{0},t_{1}-\eps)\\
                \displaystyle{\frac{t_{1}-t}{\eps}},&t\in[t_{1}-\eps,t_{1}),\\
                0,&\mbox{elsewhere,}\\
	        \end{cases}$	
\end{minipage} \hfill \begin{minipage}[t]{7cm}
                                 $\phi_{\eps}^{m}(t)=\begin{cases}
                                  \displaystyle{\frac{t-(t_{m}-\eps)}{\eps}},&t\in(t_{m}-\eps,t_{m}),\\ 
                                   0,&\mbox{elsewhere}.\\
	                           \end{cases}$	
                                   \end{minipage}\\

\begin{figure}[h]
\centering
\begin{tikzpicture}
	\draw[->]	(0,0) -- (3,0) node[anchor=south] {\scriptsize{$t$}};
	\draw[->]	(0,0) -- (0,2) node[anchor=east] {};
	\draw	(-0.1,1.5) node{{\scriptsize 1}};
	\draw	(0.5,0) node[anchor=north] {\scriptsize{$t_{k-1}$}};
	\draw	(2.5,0) node[anchor=north] {\scriptsize{$t_{k}$}};
	\draw (2.5,2pt) -- (2.5,-2pt) node[anchor=north] {};
	\draw (0.5,2pt) -- (0.5,-2pt) node[anchor=north] {};
	\draw[thick] (0,0) -- (0.5,0);
	\draw[thick] (0.5,1.5) -- (2.5,1.5);
	\draw[thick] (2.5,0) -- (3,0);
	\draw[dotted] (0,1.5) -- (0.5,1.5);
	\draw[dotted] (0.5,0) -- (0.5,1.5);
	\draw[dotted] (2.5,0) -- (2.5,1.5);
	\draw (1.6,1.9) node {\small{$\varphi^{k}$}}; 
\end{tikzpicture}
\caption{Test functions for the forward Euler method.}
\label{TestFE}
\end{figure}
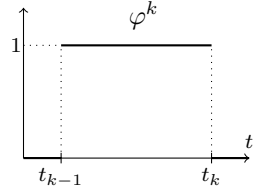

\noindent and the piecewise constant test functions are defined as (see Figure \ref{TestFE})
$$\varphi^{k}(t)=\begin{cases}
                                  1,&\;\;t\in I_{k},\\ 
                                   0,&\mbox{elsewhere},\\
	                           \end{cases}\;\forall k=1,\ldots,m.$$
	                          
We express the solution of (\ref{primaldis}) employing the trial functions defined in (\ref{trialtestFE}) as 
\begin{equation}\label{uhk} 
u_{\tau h}(t)=u_{0}\phi_{\eps}^{0}(t)+\sum_{k=1}^{m}u_{h}^{k}\phi_{\eps}^{k}(t),
\end{equation}
we approximate the initial condition $u_{0}$ by the $L^{2}$-projection in space into $V_{h}^{0}\subset V$, so we define $u_{h}^{0}$ as the solution of
$$(v_{h}^{0},u_{h}^{0})=(v_{h}^{0},u_{0}),\;\;\forall v_{h}^{0}\in V_{h}^{0}.$$ 
We commit a slight abuse of notation by omitting the constants in (\ref{uhk}) because we can express each function $u_{h}^{k}$ as a linear combination of basis functions in $V_{h}^{k},\;\forall k=1,\ldots,m.$

\begin{theorem}\label{theoFE}
(Forward Euler method) Selecting trial and test discrete subspaces as in (\ref{trialtestFE}) and taking the limit $\eps\rightarrow 0$, problem (\ref{primaldis}) leads to the following scheme 
\begin{equation}\label{primalfulldisFE}
\displaystyle{ \left\|
\begin{split}
\mbox{Find}\;u_{h}^{k}&\in V_{h}^{k},\;\forall k=1,\ldots,m,\;\mbox{such that}\\
\\
\left(v_{h}^{k},u_{h}^{k}\right)&=\left(v_{h}^{k},u_{h}^{k-1}\right)-\tau_{k}\left(\nabla v_{h}^{k},\nabla u^{k-1}_{h}\right)+\int_{I_{k}}\left<v_{h}^{k},f\right>dt,\;\forall v_{h}^{k}\in V_{h}^{k},\\
\\
\left(v_{h}^{0},u_{h}^{0}\right)&=\left(v_{h}^{0},u_{0}\right),\;\forall v_{h}^{0}\in V_{h}^{0}.
\end{split}
\right.} 
\end{equation}
which is an explicit method that is a variant of the Forward Euler method in time. 

\begin{proof}
Substituting (\ref{uhk}) into (\ref{primaldis}), we have for each test function $v_{\tau h}(t)=v_{h}^k\varphi^k(t),\forall k=1,\ldots,m$
$$B(v_{h}^{k}\varphi^{k},u_{\tau h})=B^{t}(v_{h}^{k}\varphi^{k},u_{\tau h})+B^{g}(v_{h}^{k}\varphi^{k},u_{\tau h})=\int_{I}\left<v_{h}^{k}\varphi^{k},f\right>dt.$$

As $\varphi^{k}(t)$ has local support at $I_{k}$ and its value is equal to $1$ inside the interval, we have
\begin{equation}\label{rest}
B_{|_{I_{k}}}^{t}(v_{h}^{k},u_{\tau h})+B_{|_{I_{k}}}^{g}(v_{h}^{k},u_{\tau h})=\int_{I_{k}}\left<v_{h}^{k},f\right>dt,
\end{equation}
where $B_{|_{I_{k}}}^{t}(\cdot,\cdot)$ and $B_{|_{I_{k}}}^{g}(\cdot,\cdot)$ are the restrictions of operators $B^{t}(\cdot,\cdot)$ and $B^{g}(\cdot,\cdot)$ to interval $I_{k}$, respectively. 

Figure \ref{elementFE} shows the trial shape functions over the time interval $I_{k}$, so the two left hand side terms of (\ref{rest}) become
\begin{figure}[h]
\centering
\begin{tikzpicture}
	\draw[->]	(0,0) -- (3,0) node[anchor=south] {\scriptsize{$t$}};
	\draw[->]	(0,0) -- (0,2.5) node[anchor=east] {};
	\draw	(-0.1,2) node{{\scriptsize 1}};
	\draw	(0.5,0) node[anchor=north] {\scriptsize{$t_{k-1}$}};
	\draw	(2.5,0) node[anchor=north] {\scriptsize{$t_{k}$}};
	\draw (2.5,2pt) -- (2.5,-2pt) node[anchor=north] {};
	\draw (0.5,2pt) -- (0.5,-2pt) node[anchor=north] {};
	\draw[thick] (0.5,2) -- (2,2);
	\draw[thick] (2,2) -- (2.5,0);
	\draw[thick] (2,0) -- (2.5,2);
	\draw[dotted] (2.5,0) -- (2.5,2);
	\draw[dotted] (2,2) -- (2.5,2);
	\draw[dotted] (0.5,0) -- (0.5,2);
	\draw[dotted] (0,2) -- (0.5,2);
	\draw[thick,red] (2,0) -- (2.5,0);
	\draw (2.25,0.2) node{{\textcolor{red}{\scriptsize{$\eps$}}}};
	\draw (1.3,2.5) node {\small{$\phi_{\eps}^{k-1}$}};
	\draw (2.75,2.5) node {\small{$\phi_{\eps}^{k}$}};
\end{tikzpicture}
\caption{Trial shape functions inside each time interval $I_{k}$.}
\label{elementFE}
\end{figure}
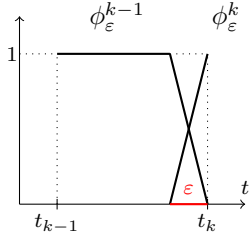
$$B_{|_{I_{k}}}^{t}(v_{h}^{k},u_{\tau h})=\int_{I_{k}}\left(v_{h}^{k},u_{h}^{k-1}\frac{\partial}{\partial t}\phi_{\eps}^{k-1}+u_{h}^{k}\frac{\partial}{\partial t}\phi_{\eps}^{k}\right)dt,$$
\begin{equation}\nonumber
\begin{split}
B_{|_{I_{k}}}^{g}(v_{h}^{k},u_{\tau h})=\int_{I_{k}}\left(\nabla v_{h}^{k},\nabla u_{h}^{k-1}\phi_{\eps}^{k-1}+\nabla u_{h}^{k}\phi_{\eps}^{k}\right)dt,
\end{split}
\end{equation}

We now calculate the integrals in each subinterval of $I_{k}=I_{k}^{1}\cup I_{k}^{\eps}$ where
$$I_{k}^{1}:=(t_{k-1},t_{k}-\eps),\;\mbox\; I_{k}^{\eps}:=(t_{k}-\eps,t_{k}),$$
that is, we separate each interval into the section where the trial functions are constant and where they linearly vary between zero and one. Then, we obtain
\begin{equation}\nonumber
\begin{split}
B_{|_{I_{k}^{1}}}^{t}(v_{h}^{k},u_{\tau h})&=\int_{t_{k-1}}^{t_{k}-\eps}\left(v_{h}^{k},u_{h}^{k-1}\cdot0\right)dt=0,\\
B_{|_{I_{k}^{\eps}}}^{t}(v_{h}^{k},u_{\tau h})&=\int_{t_{k}-\eps}^{t_{k}}\left(v_{h}^{k},u_{h}^{k-1}\left(\frac{-1}{\eps}\right)+u_{h}^{k}\frac{1}{\eps}\right)dt=(v_{h}^{k},u_{h}^{k}-u_{h}^{k-1}),\\
B_{|_{I_{k}^{1}}}^{g}(v_{h}^{k},u_{\tau h})&=\int_{t_{k-1}}^{t_{k}-\eps}(\nabla v_{h}^{k},\nabla u_{h}^{k-1})dt=(\tau_{k}-\eps)(\nabla v_{h}^{k},\nabla u_{h}^{k-1}),\\
B_{|_{I_{k}^{\eps}}}^{g}(v_{h}^{k},u_{\tau h})&=\int_{t_{k}-\eps}^{t_{k}}\left(\nabla v_{h}^{k},\nabla u_{h}^{k-1}\left(\frac{t_{k}-t}{\eps}\right)+\nabla u_{h}^{k}\left(\frac{t-(t_{k}-\eps)}{\eps}\right)\right)dt=\\
&=\frac{\eps}{2}(\nabla v_{h}^{k},\nabla u_{h}^{k-1}+\nabla u_{h}^{k}).
\end{split}
\end{equation}

Thus, local-in-time problem (\ref{rest}) becomes
\begin{equation}\nonumber
\begin{split}
\left(v_{h}^{k},u_{h}^{k}-u_{h}^{k-1}\right)&+\frac{\eps}{2}\left(\nabla v_{h}^{k},\nabla u_{h}^{k}\right)+\left(\tau_{k}-\frac{\eps}{2}\right)\left(\nabla v_{h}^{k},\nabla u_{h}^{k-1}\right)=\int_{I_{k}}\left<v_{h}^{k},f\right>dt,
\end{split}
\end{equation}
and taking the limit as $\eps\rightarrow 0$ we have
\begin{equation}\nonumber
\begin{split}
\left(v_{h}^{k},u_{h}^{k}-u_{h}^{k-1}\right)&+\tau_{k}\left(\nabla v_{h}^{k},\nabla u_{h}^{k-1}\right)=\int_{I_{k}}\left<v_{h}^{k},f\right>dt,
\end{split}
\end{equation}
$\forall k=1,\ldots,m.$
\end{proof}
\end{theorem}

Scheme (\ref{primalfulldisFE}) is the Forward Euler method in time except for the source term. A standard difference between variational forms and difference methods is that variational forms include an integral measure rather than a pointwise sample of the forcing terms. In space, we can then employ the Spectral Element Method, which leads to a diagonal mass matrix for arbitrary dimensional problems using arbitrary geometrical mappings \cite{pozrikidis2005introduction}.

\begin{remark}
To obtain an expression whose form is identical to the classical Forward Euler method, we can project the source term in the trial space as 
$$f(\mathbf{x},t)=\sum_{k=0}^{m}f^{k}(\mathbf{x})\phi_{\eps}^{k}(t),$$
where we identify $f^{k}(\mathbf{x})$ with $f(\mathbf{x},t_{k})$. Then, the source term in (\ref{primalfulldisFE}) becomes (see Figure \ref{elementFE})
$$\int_{I_{k}}\left<v_{h}^{k},f\right>dt=\int_{I_{k}}\left<v_{h}^{k},f^{k-1}\phi_{\eps}^{k-1}+f^{k}\phi_{\eps}^{k}\right>dt=\left(\tau_{k}-\frac{\eps}{2}\right)\left<v_{h}^{k},f^{k-1}\right>+\frac{\eps}{2}\left<v_{h}^{k},f^{k}\right>,$$
and taking the limit $\eps\to0$, we obtain $\forall k=1,\ldots,m$
$$\int_{I_{k}}\left<v_{h}^{k},f\right>dt=\tau_{k}\left<v_{h}^{k},f^{k-1}\right>.$$
\end{remark}

\begin{remark}
The key point in the construction of (\ref{primalfulldisFE}) is that the spatial discrete spaces in the test space $\mathcal{V}_{\tau h}$ are displaced in time with respect to the trial space $\mathcal{U}_{\tau h}^{\eps}$, which leads to a Petrov-Galerkin method. Figure \ref{spaces} illustrates this displacement of the spaces. 

\begin{figure}[h]
\centering
\begin{tikzpicture}
\draw[dotted][-] (-1,0) -- (0,0) ; 
\draw[thick][-] (0,0) -- (10,0) ; 
\draw[dotted][-] (10,0) -- (11,0) ; 
		\foreach \x in {0,2,4,8,10}
     		\draw[thick][-]  (\x,-0.1) -- (\x,0.1) node[anchor=north] {};
		\draw (0,-0.35) node{$t_{0}$};
		\draw (2,-0.35) node{$t_{1}$};
		\draw (4,-0.35) node{$t_{2}$};
		\draw (6,-0.35) node{$\ldots$};
		\draw (8,-0.35) node{$t_{m-1}$};
		\draw (10,-0.35) node{$t_{m}$};
		\draw (1,0.5) node{$V_{h}^{0}$};
		\draw (3,0.5) node{$V_{h}^{1}$};
		\draw (9,0.5) node{$V_{h}^{m-1}$};
		\draw (10.5,0.5) node{$V_{h}^{m}$};
\end{tikzpicture}

\vspace{0.5cm}

\begin{tikzpicture}
\draw[dotted][-] (-1,0) -- (0,0) ; 
\draw[thick][-] (0,0) -- (10,0) ; 
\draw[dotted][-] (10,0) -- (11,0) ; 
		\foreach \x in {0,2,4,8,10}
     		\draw[thick][-]  (\x,-0.1) -- (\x,0.1) node[anchor=north] {};
		\draw (0,-0.35) node{$t_{0}$};
		\draw (2,-0.35) node{$t_{1}$};
		\draw (4,-0.35) node{$t_{2}$};
		\draw (6,-0.35) node{$\ldots$};
		\draw (8,-0.35) node{$t_{m-1}$};
		\draw (10,-0.35) node{$t_{m}$};
		\draw (-0.5,0.5) node{$V_{h}^{0}$};
		\draw (1,0.5) node{$V_{h}^{1}$};
		\draw (3,0.5) node{$V_{h}^{2}$};
		\draw (9,0.5) node{$V_{h}^{m}$};
\end{tikzpicture}
\caption{Illustration of the displacement in time of the trial (top) and test (bottom) discrete spaces.}
\label{spaces}
\end{figure}
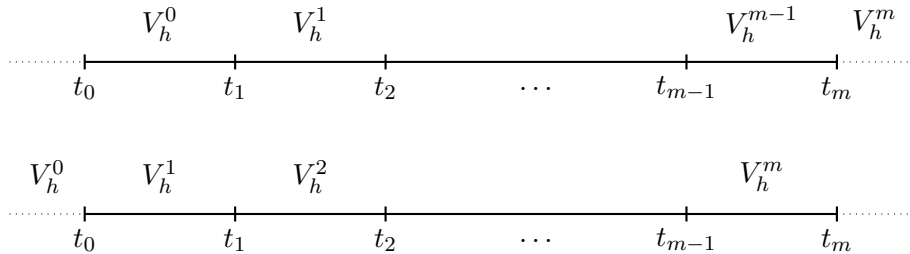
\end{remark}

\subsection{Two-stage Runge-Kutta methods}\label{S3.2}
In this section, following the same structure as in Section \ref{S3.1}, we build the trial and test spaces to obtain equivalent methods to some two-stage second-order Runge-Kutta methods. 

We introduce two test functions per temporal element $I_{k}$
\begin{equation}\label{testRK2}
v_{h,1}\varphi^{k}_{1}(t),\;v_{h,2}\varphi^{k}_{2}(t),
\end{equation}
and two trial functions associated to each element, where one of them has support over two consecutive elements (see Figure \ref{shapeRK2})
\begin{equation}\label{trialRK2}
u_{\tau h}(t)_{|_{I_{k}}}=u_{h,1}^{k-1}\phi_{\eps,1}^{k-1}(t)+u_{h,2}^{k-1}\phi_{\eps,2}^{k-1}(t)+u_{h,1}^{k}\phi_{\eps,1}^{k}(t),
\end{equation}
where the test functions could be discontinuous in time and the trial functions are globally continuous. As before, $\phi_{\eps,1}^{k}(t)$ only has support in $I_{k}^{\eps}$ and $\phi_{\eps,1}^{k-1}(t)$ and $\phi_{\eps,2}^{k-1}(t)$ are continuous and piecewise polynomials defined in $I_{k}=I_{k}^{1}\cup I_{k}^{\eps}$. We assume that the trial functions are bounded, and additionally, that each function is associated with a coefficient (as in Figure \ref{shapeRK2}). Thus, we impose
\begin{equation}\label{cond1}
\displaystyle{\left\{
\begin{split}
&\phi_{\eps,2}^{k-1}(t_{k-1})=0,\\
&\phi_{\eps,1}^{k}(t_{k}-\eps)=0,\\
&\phi_{\eps,1}^{k-1}(t_{k})=\phi_{\eps,2}^{k-1}(t_{k})=0,\\
&\phi_{\eps,1}^{k}(t_{k})=\phi_{\eps,1}^{k-1}(t_{k-1}).
\end{split}
\right.}
\end{equation}

\definecolor{forestgreen}{rgb}{0.13, 0.55, 0.13}
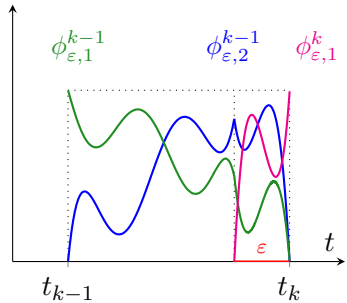
\begin{figure}[h!]
\centering
\begin{tikzpicture}
\begin{axis}[width=6cm,height=5cm,
	xmin=0,xmax=3,ymin=0,ymax=3,
	axis lines=middle,
         xlabel = $t$,
         xtick={\empty},
    	ytick={\empty},
        extra x ticks={0.5,2.5},
         extra x tick labels={$t_{k-1}$,$t_{k}$},
        ]
      \draw[dotted] (0.5,0) -- (0.5,2);
      \draw[dotted] (2.5,0) -- (2.5,2);
      \draw[dotted] (0.5,2) -- (2.5,2);
      \draw[dotted] (2,0)--(2,2);
      \draw (2.25,0.15) node{{\textcolor{red}{\scriptsize{$\eps$}}}};
      \draw (0.6,2.5) node {\textcolor{forestgreen}{\small{$\phi_{\eps,1}^{k-1}$}}};
      \draw (2,2.5) node {\textcolor{blue}{\small{$\phi_{\eps,2}^{k-1}$}}};
      \draw (2.75,2.5) node {\textcolor{magenta}{\small{$\phi_{\eps,1}^{k}$}}};
      \addplot [domain=0.5:2,samples=100,color=blue,thick]{-5.30963*x^6+62.2933*x^5-262.163*x^4+527.733*x^3-543.594*x^2+273.84*x-52.4667};
      \addplot [domain=2:2.5,samples=100,color=blue,thick]{-99.5556*x^3+657.778*x^2-1445.11*x+1057.22};
      \addplot [domain=0.5:2,samples=100,color=forestgreen,thick]{-20.8593*x^6+145.541*x^5-397.63*x^4+537.037*x^3-371.511*x^2+121.422*x-12.3333};
      \addplot [domain=2:2.5,samples=1000,color=forestgreen,thick]{-92.4444*x^3+621.333*x^2-1388.22*x+1031.67};
      \addplot [domain=2:2.5,samples=100,color=magenta,thick]{142.222*x^3-967.111*x^2+2187.11*x-1643.56};
      \draw[thick,red] (2,0) -- (2.5,0);
\end{axis}
\end{tikzpicture}
\caption{Trial functions of arbitrary order inside each element $I_{k}$.}
\label{shapeRK2}
\end{figure}

Therefore, from (\ref{cond1}), $u_{h,1}^{k-1}\in V_{h}^{k-1}$ and $u_{h,1}^{k}\in V_{h}^{k}$ are the values of $u_{\tau h}(t)$ at $t_{k-1}$ and $t_{k}$, respectively, and $u_{h,2}^{k-1}\in V_{h}^{k-1}$ is an intermediate value inside the temporal interval $I_{k}$. 
     
If we substitute (\ref{trialRK2}) into (\ref{primaldis}), for each test function, we have 
\begin{equation}\nonumber
\begin{alignedat}{6}
&&&B_{|_{I_{k}}}^{t}(v_{h,j}\varphi_{j}^{k},u_{h,1}^{k}\phi_{\eps,1}^{k})&+&B_{|_{I_{k}}}^{t}(v_{h,j}\varphi_{j}^{k},u_{h,2}^{k-1}\phi_{\eps,2}^{k-1})&+&\\ 
&+&&B_{|_{I_{k}}}^{t}(v_{h,j}\varphi_{j}^{k},u_{h,1}^{k-1}\phi_{\eps,1}^{k-1})&+&B_{|_{I_{k}}}^{g}(v_{h,j}\varphi_{j}^{k},u_{h,1}^{k}\phi_{\eps,1}^{k})&+&\\ 
&+&&B_{|_{I_{k}}}^{g}(v_{h,j}\varphi_{j}^{k},u_{h,2}^{k-1}\phi_{\eps,2}^{k-1})&+&B_{|_{I_{k}}}^{g}(v_{h,j}\varphi_{j}^{k},u_{h,1}^{k-1}\phi_{\eps,1}^{k-1})&=&\\ 
&&&&=&\int_{I_{k}}\left<v_{h,j}\varphi_{j}^{k},f\right>dt,&
\end{alignedat}
\end{equation}
where $j\in\{1,2\},$ or equivalently
\begin{equation}\label{generalRK2}
\begin{alignedat}{6}
&&&\left(v_{h,j},u_{h,1}^{k}\right)\int_{I_{k}}\varphi_{j}^{k}\frac{\partial}{\partial t}\phi_{\eps,1}^{k}\;dt&+&\left(v_{h,j},u_{h,2}^{k-1}\right)\int_{I_{k}}\varphi_{j}^{k}\frac{\partial}{\partial t}\phi_{\eps,2}^{k-1}\;dt)&+&\\
&+&&\left(v_{h,j},u_{h,1}^{k-1}\right)\int_{I_{k}}\varphi_{j}^{k}\frac{\partial}{\partial t}\phi_{\eps,1}^{k-1}\;dt&+&\left(\nabla v_{h,j},\nabla u_{h,1}^{k}\right)\int_{I_{k}}\varphi_{j}^{k}\phi_{\eps,1}^{k}\;dt&+&\\
&+&&\left(\nabla v_{h,j},\nabla u_{h,2}^{k-1}\right)\int_{I_{k}}\varphi_{j}^{k}\phi_{\eps,2}^{k-1}\;dt&+&\left(\nabla v_{h,j},\nabla u_{h,1}^{k-1}\right)\int_{I_{k}}\varphi_{j}^{k}\phi_{\eps,1}^{k-1}\;dt&=&\\
&&&&=&\int_{I_{k}}\left<v_{h,j}\varphi_{j}^{k},f\right>dt.&
\end{alignedat}
\end{equation}

We build the trial and test functions to guarantee the satisfaction of some design conditions as $\eps\rightarrow0$. We need the following orthogonality conditions in order to obtain an explicit method: 
\begin{subequations}\label{C1}
\begin{empheq}[left=\empheqlbrace]{align} 
&\lim_{\eps\to0}\int_{I_{k}}\varphi_{1}^{k}\phi_{\eps,1}^{k}\;dt=0,\label{C1.a}\\ 
&\lim_{\eps\to0}\int_{I_{k}}\varphi_{2}^{k}\frac{\partial}{\partial t}\phi_{\eps,1}^{k}\;dt=0,\label{C1.b}\\ 
&\lim_{\eps\to0}\int_{I_{k}}\varphi_{2}^{k}\phi_{\eps,1}^{k}\;dt=0,\label{C1.c}\\ 
&\lim_{\eps\to0}\int_{I_{k}}\varphi_{2}^{k}\phi_{\eps,2}^{k-1}\;dt=0.\label{C1.d}
\end{empheq}
\end{subequations}

To obtain a Runge-Kutta method we need to impose further conditions on the system. Indeed, the general expression of the two-stage second order Runge-Kutta method we want to obtain is
\begin{equation}\label{generalRK2alpha}
\begin{alignedat}{6}
&(v_{h,1},u_{h,1}^{k})&-&(v_{h,1},u_{h,1}^{k-1})+\frac{\tau_{k}}{2\alpha}(\nabla v_{h,1},\nabla u_{h,2}^{k-1})&+&\\ 
&&+&\hspace{0.3cm}\left(1-\frac{1}{2\alpha}\right)\tau_{k}(\nabla v_{h,1},\nabla u_{h,1}^{k-1})&=&\int_{I_{k}}\left<v_{h,1}\varphi_{1}^{k},f\right>dt,\\
&(v_{h,2},u_{h,2}^{k-1})&-&(v_{h,2},u_{h,1}^{k-1})+\alpha\tau_{k}(\nabla v_{h,2},\nabla u_{h,1}^{k-1})&=&\int_{I_{k}}\left<v_{h,2}\varphi_{2}^{k},f\right>dt,
\end{alignedat}
\end{equation}

where $\alpha\in\R-\{0\}$. In order to obtain (\ref{generalRK2alpha}) from (\ref{generalRK2}), in addition to the orthogonality conditions (14), we need to impose also the following conditions:
\begin{subequations}\label{C2}
\begin{empheq}[left=\empheqlbrace]{align} 
&\lim_{\eps\to0}\int_{I_{k}}\varphi_{1}^{k}\frac{\partial}{\partial t}\phi_{\eps,1}^{k}\;dt=1,\label{C2.a}\\
&\lim_{\eps\to0}\int_{I_{k}}\varphi_{1}^{k}\frac{\partial}{\partial t}\phi_{\eps,2}^{k-1}\;dt=0,\label{C2.b}\\
&\lim_{\eps\to0}\int_{I_{k}}\varphi_{1}^{k}\frac{\partial}{\partial t}\phi_{\eps,1}^{k-1}\;dt=-1,\label{C2.c}\\
&\lim_{\eps\to0}\int_{I_{k}}\varphi_{1}^{k}\phi_{\eps,2}^{k-1}\;dt=\frac{\tau_{k}}{2\alpha},\label{C2.d}\\
&\lim_{\eps\to0}\int_{I_{k}}\varphi_{1}^{k}\phi_{\eps,1}^{k-1}\;dt=\left(1-\frac{1}{2\alpha}\right)\tau_{k},\label{C2.e}\\
&\lim_{\eps\to0}\int_{I_{k}}\varphi_{2}^{k}\frac{\partial}{\partial t}\phi_{\eps,2}^{k-1}\;dt=1,\label{C2.f}\\
&\lim_{\eps\to0}\int_{I_{k}}\varphi_{2}^{k}\frac{\partial}{\partial t}\phi_{\eps,1}^{k-1}\;dt=-1,\label{C2.g}\\
&\lim_{\eps\to0}\int_{I_{k}}\varphi_{2}^{k}\phi_{\eps,1}^{k-1}\;dt=\alpha\tau_{k}.\label{C2.h}
\end{empheq}
\end{subequations}

Now, we define the limit of the shape functions as
\begin{equation}\label{limitfun}
\phi_{1}^{k-1}(t):=\lim_{\eps\to0}\phi_{\eps,1}^{k-1}(t),\;\;\phi_{2}^{k-1}(t):=\lim_{\eps\to0}\phi_{\eps,2}^{k-1}(t),\;\;\phi_{1}^{k}(t):=\lim_{\eps\to0}\phi_{\eps,1}^{k}(t),
\end{equation}
that could be discontinuous at both endpoints of $I_{k}$. We also define the jump of a function at instant $t_{k}$ as $[\phi]_{k}:=\phi(t_{k}^{+})-\phi(t_{k}^{-})$, where
$$\phi(t_{k}^{+}):=\lim_{s\longrightarrow 0^{+}}\phi(t_{k}+s),\;\;\phi(t_{k}^{-}):=\lim_{s\longrightarrow 0^{+}}\phi(t_{k}-s).$$

In the following theorem, we express the conditions (\ref{C1}) and (\ref{C2}) in the limit when $\eps\to0$. 

\begin{theorem}
When $\eps\to0$, conditions (\ref{C1}) and (\ref{C2}) are equivalent to 
\begin{subequations}\label{Cond}
\begin{alignat}{3}
&\;\;\varphi_{1}^{k}(t_{k}^{-})[\phi_{1}^{k}]_{k}=1,&&\hspace{0.4cm}\varphi_{2}^{k}(t_{k}^{-})[\phi_{1}^{k}]_{k}=0,\label{Cond.a}\\
&\int_{I_{k}}\varphi_{1}^{k}\frac{\partial}{\partial t}\phi_{1}^{k-1}\;dt+\varphi_{1}^{k}(t_{k}^{-})[\phi_{1}^{k-1}]_{k}=-1,&&\;\;\int_{I_{k}}\varphi_{1}^{k}\frac{\partial}{\partial t}\phi_{2}^{k-1}\;dt+\varphi_{1}^{k}(t_{k}^{-})[\phi_{2}^{k-1}]_{k}=0,\notag\\
&\int_{I_{k}}\varphi_{2}^{k}\frac{\partial}{\partial t}\phi_{1}^{k-1}\;dt+\varphi_{2}^{k}(t_{k}^{-})[\phi_{1}^{k-1}]_{k}=-1,&&\;\;\int_{I_{k}}\varphi_{2}^{k}\frac{\partial}{\partial t}\phi_{2}^{k-1}\;dt+\varphi_{2}^{k}(t_{k}^{-})[\phi_{2}^{k-1}]_{k}=1,\label{Cond.b}\\
&\int_{I_{k}}\varphi_{1}^{k}\phi_{1}^{k-1}\;dt=\left(1-\frac{1}{2\alpha}\right)\tau_{k},&&\;\;\int_{I_{k}}\varphi_{1}^{k}\phi_{2}^{k-1}\;dt=\frac{\tau_{k}}{2\alpha}, \notag \\
&\int_{I_{k}}\varphi_{2}^{k}\phi_{1}^{k-1}\;dt=\alpha\tau_{k},&&\;\;\int_{I_{k}}\varphi_{2}^{k}\phi_{2}^{k-1}\;dt=0. \rule{5cm}{0pt}\raisetag{36pt}\label{Cond.c}
\end{alignat}
\end{subequations}

\begin{proof}
From (\ref{C2.b}), (\ref{C2.c}), (\ref{C2.f}) and (\ref{C2.g}), we obtain
$$\lim_{\eps\to0}\int_{I_{k}}\varphi_{j}^{k}\frac{\partial}{\partial t}\phi_{\eps,i}^{k-1}\;dt=\lim_{\eps\to0}\int_{t_{k-1}}^{t_{k}-\eps}\varphi_{j}^{k}\frac{\partial}{\partial t}\phi_{\eps,i}^{k-1}\;dt+\lim_{\eps\to0}\int_{t_{k}-\eps}^{t_{k}}\varphi_{j}^{k}\frac{\partial}{\partial t}\phi_{\eps,i}^{k-1}\;dt,$$
where $i,j\in\{1,2\}$. We construct $\phi_{j}^{k}$ as continuous functions in such a way that $\displaystyle{\frac{\partial}{\partial t}\phi_{\eps,i}^{k-1}}$ do not change sign in $I_{k}^{\eps}$, so applying the first mean value theorem \cite{courant2011differential}, the second integral in the above expression becomes
$$\lim_{\eps\to0}\int_{t_{k}-\eps}^{t_{k}}\varphi_{j}^{k}\frac{\partial}{\partial t}\phi_{\eps,i}^{k-1}\;dt=\varphi_{j}^{k}(\xi)\lim_{\eps\to0}\left(\phi_{\eps,i}^{k-1}(t_{k})-\phi_{\eps,i}^{k-1}(t_{k}-\eps)\;\right),$$
were $\xi\in I_{k}^{\eps}$. Finally, taking the limit $\eps\to0$, we obtain (\ref{Cond.b})
$$\lim_{\eps\to0}\int_{I_{k}}\varphi_{j}^{k}\frac{\partial}{\partial t}\phi_{\eps,i}^{k-1}\;dt=\int_{t_{k-1}}^{t_{k}}\varphi_{j}^{k}\frac{\partial}{\partial t}\phi_{i}^{k-1}\;dt+\varphi_{j}^{k}(t_{k}^{-})[\phi_{i}^{k-1}]_{k}.$$
In a similar way, from (\ref{C1.b}) and (\ref{C2.a}), we calculate (\ref{Cond.a})
$$\lim_{\eps\to0}\int_{I_{k}}\varphi_{j}^{k}\frac{\partial}{\partial t}\phi_{\eps,1}^{k}\;dt=\lim_{\eps\to0}\int_{t_{k}-\eps}^{t_{k}}\varphi_{j}^{k}\frac{\partial}{\partial t}\phi_{\eps,1}^{k}\;dt=\varphi_{j}^{k}(t_{k}^{-})[\phi_{1}^{k}]_{k}.$$
for $j\in\{1,2\}.$

Alternatively, from (\ref{C1.d}), (\ref{C2.d}), (\ref{C2.e}) and (\ref{C2.h})
$$\lim_{\eps\to0}\int_{I_{k}}\varphi_{j}^{k}\phi_{\eps,i}^{k-1}\;dt=\lim_{\eps\to0}\int_{t_{k-1}}^{t_{k}-\eps}\varphi_{j}^{k}\phi_{\eps,i}^{k-1}\;dt+\lim_{\eps\to0}\int_{t_{k}-\eps}^{t_{k}}\varphi_{j}^{k}\phi_{\eps,i}^{k-1}\;dt,$$
for $i,j\in\{1,2\}$, but both $\varphi^{k}$ and $\phi_{\eps,i}^{k-1}$ are bounded functions over $I_{k}$, so the limit when $\eps\to0$ of the integral over $I_{k}^{\eps}$ is $0$ and we have (\ref{Cond.c})
$$\lim_{\eps\to0}\int_{I_{k}}\varphi_{j}^{k}\phi_{\eps,i}^{k-1}\;dt=\lim_{\eps\to0}\int_{t_{k-1}}^{t_{k}-\eps}\varphi_{j}^{k}\phi_{\eps,i}^{k-1}\;dt=\int_{I_{k}}\varphi_{j}^{k}\phi_{i}^{k-1}\;dt.$$
Similarly, (\ref{C1.a}) and (\ref{C1.c}) become
$$\lim_{\eps\to0}\int_{I_{k}}\varphi_{j}^{k}\phi_{\eps,1}^{k}\;dt=\lim_{\eps\to0}\int_{t_{k}-\eps}^{t_{k}}\varphi_{j}^{k}\phi_{\eps,1}^{k}\;dt=0.$$
\end{proof}
\end{theorem}

From (\ref{cond1}), functions (\ref{limitfun}) satisfy
\begin{equation}\label{cond2}
\displaystyle{\left\{
\begin{split}
&\phi_{2}^{k-1}(t_{k-1}^{+})=0,\\
&\phi_{1}^{k}(t_{k}^{-})=0,\\
&\phi_{1}^{k-1}(t_{k}^{+})=\phi_{2}^{k-1}(t_{k}^{+})=0,\\
&\phi_{1}^{k}(t_{k}^{+})=\phi_{1}^{k-1}(t_{k-1}^{+}),\\
\end{split}
\right.}
\end{equation}
then conditions (\ref{Cond}) become
\begin{subequations}\label{Condlim}
\begin{alignat}{3}
&\;\;\varphi_{1}^{k}(t_{k}^{-})\phi_{1}^{k-1}(t_{k-1}^{+})=1,&&\hspace{0.4cm}\varphi_{2}^{k}(t_{k}^{-})\phi_{1}^{k-1}(t_{k-1}^{+})=0,\label{Condlim.a}\\
&\int_{I_{k}}\varphi_{1}^{k}\frac{\partial}{\partial t}\phi_{1}^{k-1}\;dt-\varphi_{1}^{k}(t_{k}^{-})\phi_{1}^{k-1}(t_{k}^{-})=-1,&&\;\;\int_{I_{k}}\varphi_{1}^{k}\frac{\partial}{\partial t}\phi_{2}^{k-1}\;dt-\varphi_{1}^{k}(t_{k}^{-})\phi_{2}^{k-1}(t_{k}^{-})=0,\notag\\
&\int_{I_{k}}\varphi_{2}^{k}\frac{\partial}{\partial t}\phi_{1}^{k-1}\;dt-\varphi_{2}^{k}(t_{k}^{-})\phi_{1}^{k-1}(t_{k}^{-})=-1,&&\;\;\int_{I_{k}}\varphi_{2}^{k}\frac{\partial}{\partial t}\phi_{2}^{k-1}\;dt-\varphi_{2}^{k}(t_{k}^{-})\phi_{2}^{k-1}(t_{k}^{-})=1,\label{Condlim.b}\\
&\int_{I_{k}}\varphi_{1}^{k}\phi_{1}^{k-1}\;dt=\left(1-\frac{1}{2\alpha}\right)\tau_{k},&&\;\;\int_{I_{k}}\varphi_{1}^{k}\phi_{2}^{k-1}\;dt=\frac{\tau_{k}}{2\alpha},\notag\\
&\int_{I_{k}}\varphi_{2}^{k}\phi_{1}^{k-1}\;dt=\alpha\tau_{k},&&\;\;\int_{I_{k}}\varphi_{2}^{k}\phi_{2}^{k-1}\;dt=0.\rule{5cm}{0pt}\raisetag{36pt}\label{Condlim.c}
\end{alignat}
\end{subequations}

From (\ref{Condlim.a}), we obtain $\varphi_{2}^{k}(t_{k}^{-})=0$, and we assume that $\varphi_{1}^{k}(t_{k}^{-})=1$, so the conditions that the function must satisfy together with (\ref{cond2}), become
\begin{subequations}\label{Condfinal}
\begin{alignat}{3}
&\;\;\varphi_{1}^{k}(t_{k}^{-})=1,&&\hspace{0.4cm}\varphi_{2}^{k}(t_{k}^{-})=0,\notag \\
&\;\;\phi_{1}^{k-1}(t_{k-1}^{+})=1,&&\hspace{0.4cm}\phi_{2}^{k-1}(t_{k-1}^{+})=0,\label{Condfinal.a} \\
&\int_{I_{k}}\varphi_{1}^{k}\frac{\partial}{\partial t}\phi_{1}^{k-1}\;dt-\phi_{1}^{k-1}(t_{k}^{-})=-1,&&\;\;\int_{I_{k}}\varphi_{1}^{k}\frac{\partial}{\partial t}\phi_{2}^{k-1}\;dt-\phi_{2}^{k-1}(t_{k}^{-})=0,\notag\\
&\int_{I_{k}}\varphi_{2}^{k}\frac{\partial}{\partial t}\phi_{1}^{k-1}\;dt=-1,&&\;\;\int_{I_{k}}\varphi_{2}^{k}\frac{\partial}{\partial t}\phi_{2}^{k-1}\;dt=1,  \label{Condfinal.b}\\
&\int_{I_{k}}\varphi_{1}^{k}\phi_{1}^{k-1}\;dt=\left(1-\frac{1}{2\alpha}\right)\tau_{k},&&\;\;\int_{I_{k}}\varphi_{1}^{k}\phi_{2}^{k-1}\;dt=\frac{\tau_{k}}{2\alpha},\notag\\
&\int_{I_{k}}\varphi_{2}^{k}\phi_{1}^{k-1}\;dt=\alpha\tau_{k},&&\;\;\int_{I_{k}}\varphi_{2}^{k}\phi_{2}^{k-1}\;dt=0.\rule{5cm}{0pt}\raisetag{36pt}\label{Condfinal.c}
\end{alignat}
\end{subequations}

For the trial and test functions which are polynomials of arbitrary order, then (\ref{Condfinal}) becomes a system of nonlinear equations. In particular, if we select linear trial and test functions, we obtain no solutions. However, if we select quadratic functions, we have a system of $12$ nonlinear equations with $12$ unknowns that has two solutions. We solve the resulting system in the master element $[0,1]$ using the MATLAB code we describe in \ref{A2} (see Table \ref{solRK2}), and we obtain different sets of trial and test functions depending on the value of $\alpha$.

\renewcommand{\arraystretch}{1.5}
\begin{table}[h!]
\centering
\begin{tabular}{|c|c|c|}\hline
&Solution 1&Solution 2\\\hline
$\phi_{1}^{k-1}(t)$&$-\frac{1}{\alpha}t+1$&$-\frac{6}{\alpha}t^{2}+\frac{3}{\alpha}t+1$\\
$\phi_{2}^{k-1}(t)$&$\frac{1}{\alpha}t$&$\frac{6}{\alpha}t^{2}-\frac{3}{\alpha}t$\\
$\varphi_{1}^{k}(t)$&$1$&$1$\\
$\varphi_{2}^{k}(t)$&$12\alpha t^{2}-18\alpha t+6\alpha$&$-2\alpha t+2\alpha$\\\hline
\end{tabular}
\caption{Trial and test functions defined over the master element $[0,1]$ that lead to the two-stage Runge-Kutta method.}
\label{solRK2}
\end{table}

\subsection*{Example 1: Explicit trapezoidal rule}
When $\alpha=1$, (\ref{generalRK2alpha}) is equivalent to the \textit{Explicit trapezoidal rule} \cite{hairer2006geometric}. Figure \ref{Heun} shows the trial and test functions of both solutions over the master element. 

\pgfmathsetmacro{\width} {5cm}
\pgfmathsetmacro{\height} {5cm}
\begin{figure}[h!]
\begin{subfigure}[b]{0.5\linewidth}
\centering
\begin{tikzpicture}
\begin{axis}[footnotesize,
		  width=\width,height=\height, xlabel=$t$,
		  xtick={0,1},ytick={0,1},
		  enlargelimits=false]
\addplot [domain=0:1,samples=100,color=blue,thick]{1-x};
\addplot [domain=0:1,samples=100,color=red,thick]{x};
\end{axis}
\end{tikzpicture}
\caption{\footnotesize Trial functions of the first solution.} 
\end{subfigure}
\begin{subfigure}[b]{0.5\linewidth}
\centering
\begin{tikzpicture}
\begin{axis}[footnotesize,
		  width=\width,height=\height, xlabel=$t$,
		  xtick={0,1},ytick={0,1,6},
		  enlargelimits=false]
\addplot [domain=0:1,samples=100,color=blue,thick]{1};
\addplot [domain=0:1,samples=100,color=red,thick]{12*x*x-18*x+6};
\end{axis}
\end{tikzpicture}
\caption{\footnotesize Test functions of the first solution.} 
\end{subfigure}
\begin{subfigure}[b]{0.5\linewidth}
\centering
\begin{tikzpicture}
\begin{axis}[footnotesize,
		  width=\width,height=\height, xlabel=$t$,
		 xtick={0,1},ytick={-2,0,1,3},
		  enlargelimits=false]
\addplot [domain=0:1,samples=100,color=blue,thick]{-6*x*x+3*x+1};
\addplot [domain=0:1,samples=100,color=red,thick]{6*x*x-3*x};
\end{axis}
\end{tikzpicture}
\caption{\footnotesize Trial functions of the second solution.}  
\end{subfigure}
\begin{subfigure}[b]{0.5\linewidth}
\centering
\begin{tikzpicture}
\begin{axis}[footnotesize,
		  width=\width,height=\height, xlabel=$t$,
		 xtick={0,1},ytick={0,1,2},
		 enlargelimits=false]
\addplot [domain=0:1,samples=100,color=blue,thick]{1};
\addplot [domain=0:1,samples=100,color=red,thick]{-2*x+2};
\end{axis}
\end{tikzpicture}
\caption{\footnotesize Test functions of the second solution.} 
\end{subfigure}
\caption{Trial and test functions over the master element $[0,1]$ when $\alpha=1$.}
\label{Heun}
\end{figure}
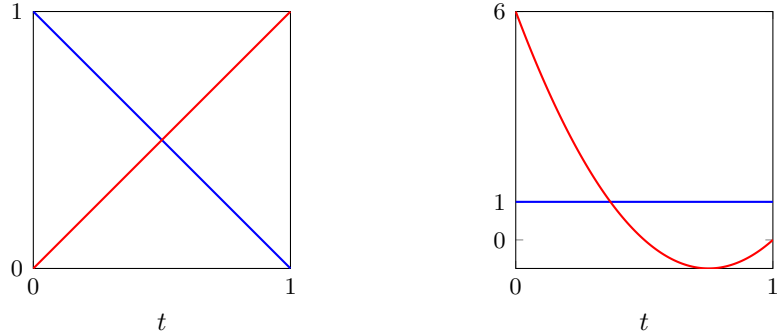
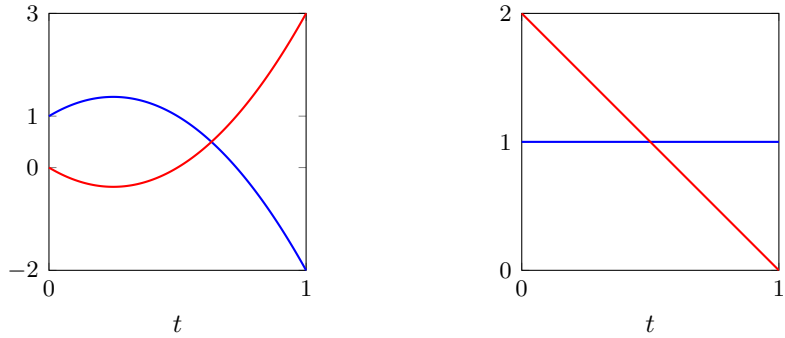

\subsection*{Example 2: Explicit midpoint rule}
When $\alpha=\frac{1}{2}$, we obtain the \textit{Explicit midpoint rule} \cite{hairer2006geometric} and in Figure \ref{Midpoint}, we can see the trial and test functions of both solutions over $[0,1]$. 

\pgfmathsetmacro{\width} {5cm}
\pgfmathsetmacro{\height} {5cm}
\begin{figure}[h!]
\begin{subfigure}[b]{0.5\linewidth}
\centering
\begin{tikzpicture}
\begin{axis}[footnotesize,
		  width=\width,height=\height, xlabel=$t$,
		  xtick={0,1},ytick={-1,0,1,2},
		  enlargelimits=false]
\addplot [domain=0:1,samples=100,color=blue,thick]{1-2*x};
\addplot [domain=0:1,samples=100,color=red,thick]{2*x};
\end{axis}
\end{tikzpicture}
\caption{\footnotesize Trial functions of the first solution.} 
\end{subfigure}
\begin{subfigure}[b]{0.5\linewidth}
\centering
\begin{tikzpicture}
\begin{axis}[footnotesize,
		  width=\width,height=\height, xlabel=$t$,
		  xtick={0,1},ytick={0,1,3},
		  enlargelimits=false]
\addplot [domain=0:1,samples=100,color=blue,thick]{1};
\addplot [domain=0:1,samples=100,color=red,thick]{6*x*x-9*x+3};
\end{axis}
\end{tikzpicture}
\caption{\footnotesize Test functions of the first solution.} 
\end{subfigure}
\begin{subfigure}[b]{0.5\linewidth}
\centering
\begin{tikzpicture}
\begin{axis}[footnotesize,
		  width=\width,height=\height, xlabel=$t$,
		 xtick={0,1},ytick={-4,0,1,6},
		  enlargelimits=false]
\addplot [domain=0:1,samples=100,color=blue,thick]{-12*x*x+6*x+1};
\addplot [domain=0:1,samples=100,color=red,thick]{12*x*x-6*x};
\end{axis}
\end{tikzpicture}
\caption{\footnotesize Trial functions of the second solution.}  
\end{subfigure}
\begin{subfigure}[b]{0.5\linewidth}
\centering
\begin{tikzpicture}
\begin{axis}[footnotesize,
		  width=\width,height=\height, xlabel=$t$,
		 xtick={0,1},ytick={0,1,2},
		 enlargelimits=false]
\addplot [domain=0:1,samples=100,color=blue,thick]{1};
\addplot [domain=0:1,samples=100,color=black]{1.1};
\addplot [domain=0:1,samples=100,color=red,thick]{-x+1};
\end{axis}
\end{tikzpicture}
\caption{\footnotesize Test functions of the second solution.} 
\end{subfigure}
\caption{Trial and test functions over the master element $[0,1]$ when $\alpha=\frac{1}{2}$.}
\label{Midpoint}
\end{figure}
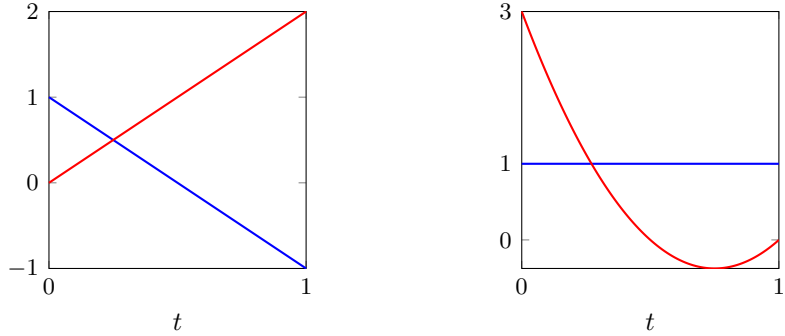
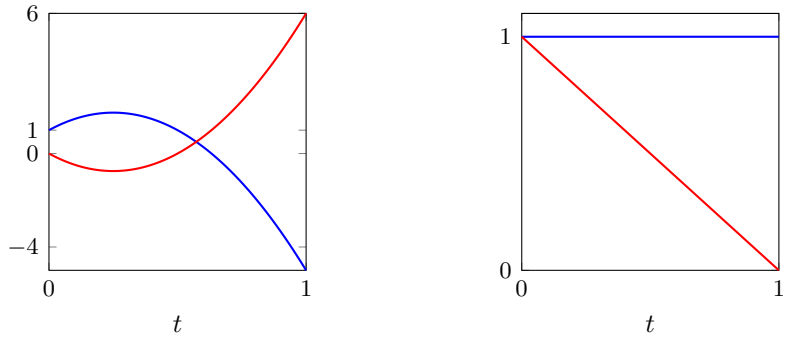

\begin{remark}
In order for the discrete system to make sense and thus result in square mass matrices in (\ref{generalRK2alpha}), we need test functions as in (\ref{testRK2}) satisfying
$$v_{h,1}\in V_{h}^{k},\;\;v_{h,2}\in V_{h}^{k-1},$$
so while both test functions are polynomials in time, in space they belong to different spaces. 
\end{remark}

\begin{remark}
As in Section \ref{S3.1}, to obtain expressions whose form is identical to standard Runge-Kutta methods, we can project the source term in the trial space as 
\begin{equation}\nonumber
\begin{alignedat}{4}
&\int_{I_{k}}\left<v_{h,1}\varphi_{1}^{k},f\right>dt&&=\int_{I_{k}}\left<v_{h,1}\varphi_{1}^{k},f_{1}^{k-1}\phi_{1}^{k-1}+f_{2}^{k-1}\phi_{2}^{k-1}\right>dt=\\
&&&=\left(1-\frac{1}{2\alpha}\right)\tau_{k}\left<v_{h,1},f_{1}^{k-1}\right>+\frac{\tau_{k}}{2\alpha}\left<v_{h,1},f_{2}^{k-1}\right>,
\end{alignedat}
\end{equation}
$$\int_{I_{k}}\left<v_{h,2}\varphi_{2}^{k},f\right>dt=\int_{I_{k}}\left<v_{h,2}\varphi_{2}^{k},f_{1}^{k-1}\phi_{1}^{k-1}+f_{2}^{k-1}\phi_{2}^{k-1}\right>dt=\alpha\tau_{k}\left<v_{h,2},f_{1}^{k-1}\right>,$$
which is the general two-stage Runge-Kutta method in time. 
\end{remark}

\subsection{General s-stage Runge-Kutta Methods}\label{S3.3}
In this section, we generalize the constructions of Sections \ref{S3.1} and \ref{S3.2}, to the general $s$-stage explicit Runge-Kutta Method. 

The general method we want to obtain is of the form
\begin{equation}\label{highorderRK}
\begin{alignedat}{9}
&(v_{h,1},u_{h,1}^{k})&\hspace{0.1cm}-\hspace{0.1cm}&(v_{h,1},u_{h,1}^{k-1})&\hspace{0.1cm}+\hspace{0.1cm}&\hspace{0.1cm}\tau_{k}\sum_{i=1}^{s}b_{i}(\nabla v_{h,1},\nabla u_{h,i}^{k-1})&=\hspace{0.1cm}&\tau_{k}\sum_{i=1}^{s}b_{i}\left<v_{h,1},f_{i}^{k-1}\right>,\\
&(v_{h,2},u_{h,2}^{k-1})&\hspace{0.1cm}-\hspace{0.1cm}&(v_{h,2},u_{h,1}^{k-1})&\hspace{0.1cm}+\hspace{0.1cm}&\hspace{0.1cm}\tau_{k}a_{21}(\nabla v_{h,2},\nabla u_{h,1}^{k-1})&=\hspace{0.1cm}&\tau_{k}a_{21}\left<v_{h,2},f_{1}^{k-1}\right>,\\
&(v_{h,3},u_{h,3}^{k-1})&\hspace{0.1cm}-\hspace{0.1cm}&(v_{h,3},u_{h,1}^{k-1})&\hspace{0.1cm}+\hspace{0.1cm}&\hspace{0.1cm}\tau_{k}a_{31}(\nabla v_{h,3},\nabla u_{h,1}^{k-1})&+\hspace{0.1cm}&\tau_{k}a_{32}(\nabla v_{h,3},\nabla u_{h,2}^{k-1})=\\
&&&&\hspace{0.1cm}=\hspace{0.1cm}&\hspace{0.1cm}\tau_{k}a_{31}\left<v_{h,3},f_{1}^{k-1}\right>&+\hspace{0.1cm}&\tau_{k}a_{32}\left<v_{h,3},f_{2}^{k-1}\right>,\\
&&&&\vdots\hspace{0.2cm}&&&&\\
&(v_{h,s},u_{h,s}^{k})&\hspace{0.1cm}-\hspace{0.1cm}&(v_{h,s},u_{h,s}^{k-1})&\hspace{0.1cm}+\hspace{0.1cm}&\hspace{0.1cm}\tau_{k}\sum_{j=1}^{s-1}a_{sj}(\nabla v_{h,s},\nabla u_{h,j}^{k-1})&\hspace{0.1cm}=\hspace{0.1cm}&\tau_{k}\sum_{j=1}^{s-1}a_{sj}\left<v_{h,s},f_{j}^{k-1}\right>,
\end{alignedat}
\end{equation}
or in compact form 
\begin{equation}\label{highorderRKcompact}
\begin{alignedat}{6}
&(v_{h,1},u_{h,1}^{k})&\hspace{0.1cm}-\hspace{0.1cm}&(v_{h,1},u_{h,1}^{k-1})&\hspace{0.1cm}+\hspace{0.1cm}&\tau_{k}\sum_{i=1}^{s}b_{i}(\nabla v_{h,1},\nabla u_{h,i}^{k-1})&\hspace{0.1cm}=\hspace{0.1cm}&\tau_{k}\sum_{i=1}^{s}b_{i}\left<v_{h,1},f_{i}^{k-1}\right>,\\
&(v_{h,i},u_{h,i}^{k})&\hspace{0.1cm}-\hspace{0.1cm}&(v_{h,i},u_{h,i}^{k-1})&\hspace{0.1cm}+\hspace{0.1cm}&\tau_{k}\sum_{j=1}^{i-1}a_{ij}(\nabla v_{h,i},\nabla u_{h,j}^{k-1})&\hspace{0.1cm}=\hspace{0.1cm}&\tau_{k}\sum_{j=1}^{i-1}a_{ij}\left<v_{h,i},f_{j}^{k-1}\right>,
\end{alignedat}
\end{equation}
$\;\forall i=2,\ldots,s,$ where
\begin{equation}
f_{i}^{k-1}(\mathbf{x}):=f(\mathbf{x},t_{k-1}+c_{i}\tau_{k}),\;\forall i=1,\ldots,s.
\end{equation}
The coefficients $a_{ij}$, $b_{i}$, $c_{i}$, with $i,j\in\{1,\ldots,s\}$, are the ones corresponding to the Butcher tableau (see Table \ref{Butcher}) \cite{butcher2008numerical}. As (\ref{highorderRK}) is an explicit method, we have that
$$a_{ij}=0,\forall j\geq i.$$

\begin{table}[h!]
\centering
\begin{tabular}{c|cccc}
$c_{1}$&$a_{11}$&$a_{12}$&$\ldots$&$a_{1s}$\\
$c_{2}$&$a_{21}$&$a_{22}$&$\ldots$&$a_{2s}$\\
$\vdots$&$\vdots$&$\vdots$&$\ddots$&$\vdots$\\
$c_{s}$&$a_{s1}$&$a_{s2}$&$\ldots$&$a_{ss}$\\\hline
&$b_{1}$&$b_{2}$&\ldots&$b_{s}$\\
\end{tabular}
\caption{Butcher tableau}
\label{Butcher}
\end{table}

As before, we consider $s$ trial functions
$$\phi_{1}^{k-1}(t),\ldots,\phi_{s}^{k-1}(t),$$
and $s$ test functions per time interval
$$v_{h,1}\varphi_{1}^{k}(t),\ldots,v_{h,s}\varphi_{s}^{k}(t),$$
and to obtain an expression whose form is identical to the classical Runge-Kutta methods, we project the source term over the trial space as
$$f(\mathbf{x},t)=\sum_{k=1}^{m}\sum_{i=1}^{s}f_{i}^{k-1}(\mathbf{x})\phi_{i}^{k-1}(t).$$ 

Following the same logic as in Section \ref{S3.2}, we generalize conditions (\ref{Condfinal}) to $s$-stages as follows
\begin{subequations}\label{CondfinalRK}
\begin{alignat}{3}
&\;\;\varphi_{1}^{k}(t_{k}^{-})=1,&&\hspace{0.2cm}\varphi_{j}^{k}(t_{k}^{-})=0,\;\;\forall j=2,\ldots,s,\notag\\
&\;\;\phi_{1}^{k-1}(t_{k-1}^{+})=1,&&\hspace{0.2cm}\phi_{i}^{k-1}(t_{k-1}^{+})=0,\;\;\forall i=2,\ldots,s, \rule{5cm}{0pt}\raisetag{25pt}\label{CondfinalRK.a}\\
&\int_{I_{k}}\varphi_{1}^{k}\frac{\partial}{\partial t}\phi_{1}^{k-1}\;dt-\phi_{1}^{k-1}(t_{k}^{-})=-1,&&\int_{I_{k}}\varphi_{1}^{k}\frac{\partial}{\partial t}\phi_{i}^{k-1}\;dt-\phi_{i}^{k-1}(t_{k}^{-})=0,\;\;\forall i=2,\ldots,s, \notag\\
&\int_{I_{k}}\varphi_{j}^{k}\frac{\partial}{\partial t}\phi_{1}^{k-1}\;dt=-1,&&\int_{I_{k}}\varphi_{j}^{k}\frac{\partial}{\partial t}\phi_{i}^{k-1}\;dt=\begin{cases}1,\;\mbox{if}\;i=j\\0,\;\mbox{if}\;i\neq j\end{cases}\hspace{-0.3cm},\;\;\forall i,j=2,\ldots,s,\label{CondfinalRK.b}\\
&\int_{I_{k}}\varphi_{1}^{k}\phi_{i}^{k-1}\;dt=\tau_{k}b_{i},\;\forall i=1,\ldots,s,\notag\\
&\int_{I_{k}}\varphi_{j}^{k}\phi_{i}^{k-1}\;dt=\tau_{k}a_{ji},\;\forall i=1,\ldots,s,&&\;\forall j=2,\ldots,s.\rule{5cm}{0pt}\raisetag{36pt}\label{CondfinalRK.c}
\end{alignat}
\end{subequations}

\begin{remark}
As before, to properly define solvable discrete systems, we seek to obtain square mass matrices in (\ref{highorderRK}). Thus, we need test functions satisfying
$$v_{h,1}\in V_{h}^{k},\;\;v_{h,j}\in V_{h}^{k-1},\;\;\forall j=2,\ldots,s.$$
\end{remark}

\subsection*{Example 1: Three-stage Runge-Kutta}
We calculate the trial and test functions of the three-stage and third order Runge-Kutta method that has the Butcher tableau as in Table \ref{ButcherRK3}  \cite{butcher2008numerical}.

\begin{table}[h!]
\centering
\begin{tabular}{c|cccc}
$0$&$0$&$0$&$0$\\
$\frac{1}{2}$&$\frac{1}{2}$&$0$&$0$\\
$1$&$-1$&$2$&$0$\\ \hline
&$\frac{1}{6}$&$\frac{2}{3}$&$\frac{1}{6}$\\
\end{tabular}
\caption{Butcher tableau of the three-stage Runge-Kutta method.}
\label{ButcherRK3}
\end{table}

If we consider cubic polynomials for the trial and test functions, we obtain four possible solutions: two of them with real coefficients (see Table \ref{solrealRK3}) and an the remaining two solutions with complex conjugate coefficients (see Table \ref{solcomplexRK3}). Figure \ref{RK3} shows the trial and test functions of the real solutions over the master element $[0,1]$. In Table \ref{solcomplexRK3}, $z_{j}$ and $\bar{z}_{j}$ denote the following complex numbers and their conjugates

\begin{center}
\begin{tabular}{ll}
$z_{1}=(9-i\sqrt{66})/7,$&$z_{6}=(34+3i\sqrt{66})/7,$\\
$z_{2}=(11-2i\sqrt{66})/7,$&$z_{7}=(16+i\sqrt{66})/7,$\\
$z_{3}=(12+i\sqrt{66})/7,$&$z_{8}=(30+i\sqrt{66})/7,$\\
$z_{4}=(5+i\sqrt{66})/7,$&$z_{9}=(89+6i\sqrt{66})/7,$\\
$z_{5}=(2-i\sqrt{66})/7,$&$z_{10}=(39+2i\sqrt{66})/7,$\\
\end{tabular}
\end{center}

\renewcommand{\arraystretch}{1.5}
\begin{table}[h!]
\centering
\begin{tabular}{|c|c|c|}\hline
&Solution 1&Solution 2\\\hline
$\phi_{1}^{k-1}(t)$&$\frac{1}{2}t^{2}-2t+1$&$110t^{3}-130t^{2}+30t+1$\\
$\phi_{2}^{k-1}(t)$&$-t^{2}+2t$&$-120t^{3}+140t^{2}-32t$\\
$\phi_{3}^{k-1}(t)$&$\frac{1}{2}t^{2}$&$10t^{3}-10t^{2}+2t$\\
$\varphi_{1}^{k}(t)$&$1$&$1$\\
$\varphi_{2}^{k}(t)$&$-30t^{3}+60t^{2}-36t+6$&$\frac{3}{2}t^{2}-3t+\frac{3}{2}$\\
$\varphi_{3}^{k}(t)$&$420t^{3}-780t^{2}+408t-48$&$-6t^{2}+6t$\\\hline
\end{tabular}
\caption{Trial and test functions with real coefficients over the master element $[0,1]$ that lead to a three-stage Runge-Kutta method.}
\label{solrealRK3}
\end{table}

\renewcommand{\arraystretch}{1.5}
\begin{table}[h!]
\centering
\begin{tabular}{|c|c|c|}\hline
&Solution 3&Solution 4\\\hline
$\phi_{1}^{k-1}(t)$&$\frac{10}{3}z_{1}t^{3}-2z_{2}t^{2}-z_{3}t+1$&$\frac{10}{3}\bar{z}_{1}t^{3}-2\bar{z}_{2}t^{2}-\bar{z}_{3}t+1$\\
$\phi_{2}^{k-1}(t)$&$-\frac{20}{3}z_{1}t^{3}+4z_{2}t^{2}+2z_{4}t$&$-\frac{20}{3}\bar{z}_{1}t^{3}+4\bar{z}_{2}t^{2}+2\bar{z}_{4}t$\\
$\phi_{3}^{k-1}(t)$&$\frac{10}{3}z_{1}t^{3}-2z_{2}t^{2}+z_{5}t$&$\frac{10}{3}\bar{z}_{1}t^{3}-2\bar{z}_{2}t^{2}+\bar{z}_{5}t$\\
$\varphi_{1}^{k}(t)$&$1$&$1$\\
$\varphi_{2}^{k}(t)$&$-10\bar{z}_{1}t^{3}+6z_{6}t^{2}-9z_{7}t+z_{8}$&$-10z_{1}t^{3}+6\bar{z}_{6}t^{2}-9\bar{z}_{7}t+\bar{z}_{8}$\\
$\varphi_{3}^{k}(t)$&$40\bar{z}_{1}t^{3}-24z_{6}t^{2}+6z_{9}t-2z_{10}$&$40z_{1}t^{3}-24\bar{z}_{6}t^{2}+6\bar{z}_{9}t-2\bar{z}_{10}$\\\hline
\end{tabular}
\caption{Trial and test functions with complex coefficients over the master element $[0,1]$ that lead to a three-stage Runge-Kutta method.}
\label{solcomplexRK3}
\end{table}

\pgfmathsetmacro{\width} {5cm}
\pgfmathsetmacro{\height} {5cm}
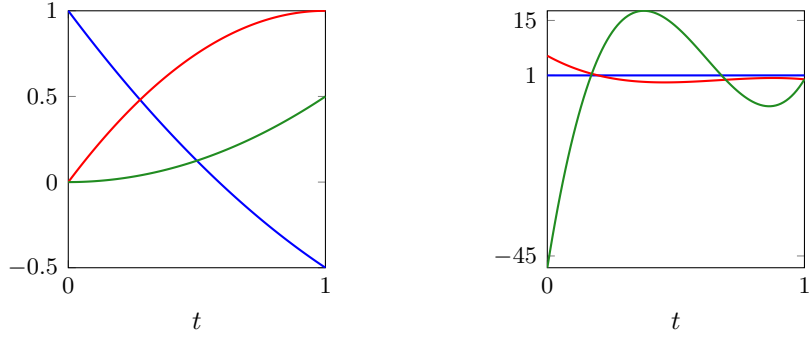
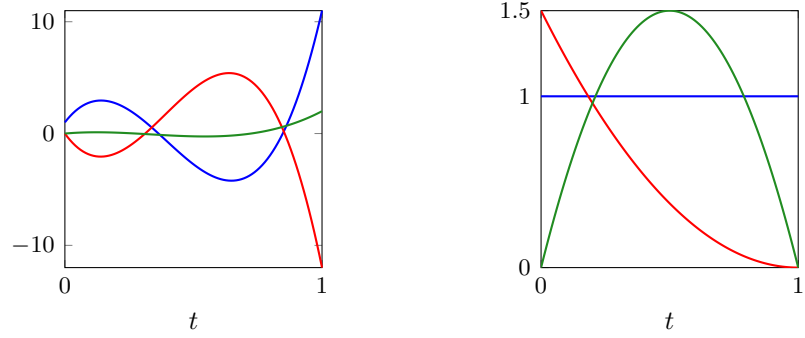
\begin{figure}[h!]
\begin{subfigure}[b]{0.5\linewidth}
\centering
\begin{tikzpicture}
\begin{axis}[footnotesize,
		  width=\width,height=\height, xlabel=$t$,
		 xtick={0,1},ytick={-0.5,0,0.5,1},
		  enlargelimits=false]
\addplot [domain=0:1,samples=100,color=blue,thick]{0.5*x*x-2*x+1};
\addplot [domain=0:1,samples=100,color=red,thick]{-x*x+2*x};
\addplot [domain=0:1,samples=100,color=forestgreen,thick]{0.5*x*x};
\end{axis}
\end{tikzpicture}
\caption{\footnotesize Trial functions of the first real solution.}  
\end{subfigure}
\begin{subfigure}[b]{0.5\linewidth}
\centering
\begin{tikzpicture}
\begin{axis}[footnotesize,
		  width=\width,height=\height, xlabel=$t$,
		 xtick={0,1},ytick={-45,1,15},
		 enlargelimits=false]
\addplot [domain=0:1,samples=100,color=blue,thick]{1};
\addplot [domain=0:1,samples=100,color=red,thick]{-30*x*x*x+60*x*x-36*x+6};
\addplot [domain=0:1,samples=100,color=forestgreen,thick]{420*x*x*x-780*x*x+408*x-48};
\end{axis}
\end{tikzpicture}
\caption{\footnotesize Test functions of the first real solution.} 
\end{subfigure}
\begin{subfigure}[b]{0.5\linewidth}
\centering
\begin{tikzpicture}
\begin{axis}[footnotesize,
		  width=\width,height=\height, xlabel=$t$,
		  xtick={0,1},ytick={-10,0,10},
		  enlargelimits=false]
\addplot [domain=0:1,samples=100,color=blue,thick]{110*x*x*x-130*x*x+30*x+1};
\addplot [domain=0:1,samples=100,color=red,thick]{-120*x*x*x+140*x*x-32*x};
\addplot [domain=0:1,samples=100,color=forestgreen,thick]{10*x*x*x-10*x*x+2*x};
\end{axis}
\end{tikzpicture}
\caption{\footnotesize Trial functions of the second real solution.} 
\end{subfigure}
\begin{subfigure}[b]{0.5\linewidth}
\centering
\begin{tikzpicture}
\begin{axis}[footnotesize,
		  width=\width,height=\height, xlabel=$t$,
		  xtick={0,1},ytick={0,1,1.5},
		  enlargelimits=false]
\addplot [domain=0:1,samples=100,color=blue,thick]{1};
\addplot [domain=0:1,samples=100,color=red,thick]{1.5*x*x-3*x+1.5};
\addplot [domain=0:1,samples=100,color=forestgreen,thick]{-6*x*x+6*x};
\end{axis}
\end{tikzpicture}
\caption{\footnotesize Test functions of the second real solution.} 
\end{subfigure}
\caption{Trial and test functions over the master element $[0,1]$ for the real solutions of the three-stage Runge-Kutta method.}
\label{RK3}
\end{figure}

\subsection*{Example 2: Four-stage Runge-Kutta}
Now, we consider the four-stage and fourth order Runge-Kutta method with the Butcher tableau as in Table \ref{ButcherRK4} \cite{butcher2008numerical}. Table \ref{solrealRK4} and Figure \ref{RK4} show two real solutions for this method.

\begin{table}[h!]
\centering
\begin{tabular}{c|cccc}
$0$&$0$&$0$&$0$&$0$\\
$\frac{1}{2}$&$\frac{1}{2}$&$0$&$0$&$0$\\
$\frac{1}{2}$&$0$&$\frac{1}{2}$&$0$&$0$\\
$1$&$0$&$0$&$1$&$0$\\ \hline
&$\frac{1}{6}$&$\frac{1}{3}$&$\frac{1}{3}$&$\frac{1}{6}$\\
\end{tabular}
\caption{Butcher tableau of the four-stage Runge-Kutta method.} 
\label{ButcherRK4}
\end{table}

\renewcommand{\arraystretch}{1.5}
\begin{table}[h!]
\centering
\begin{adjustbox}{max width=\textwidth}
\begin{tabular}{|c|c|c|}\hline
&Solution 1&Solution 2\\\hline
$\phi_{1}^{k-1}(t)$&$\frac{2}{3}t^{3}-2t+1$&$-\frac{2800}{3}t^{4}+1610t^{3}-810t^{2}+\frac{320}{3}t+1$\\
$\phi_{2}^{k-1}(t)$&$-2t^{2}+2t$&$\frac{4480}{3}t^{4}-2520t^{3}+1230t^{2}-\frac{470}{3}t$\\
$\phi_{3}^{k-1}(t)$&$-\frac{4}{3}t^{3}+2t^{2}$&$-\frac{1820}{3}t^{4}+980t^{3}-450t^{2}+\frac{160}{3}t$\\
$\phi_{4}^{k-1}(t)$&$\frac{2}{3}t^{3}$&$\frac{140}{3}t^{4}-70t^{3}+30t^{2}-\frac{10}{3}t$\\
$\varphi_{1}^{k}(t)$&$1$&$1$\\
$\varphi_{2}^{k}(t)$&$140t^{4}-350t^{3}+300t-100t+10$&$-2t^{3}+6t^{2}-6t+2$\\
$\varphi_{3}^{k}(t)$&$-910t^{4}+2170t^{3}-1725t^{2}+50t-35$&$2t^{3}-3t^{2}+1$\\
$\varphi_{4}^{k}(t)$&$4480t^{4}-10360t^{3}+7890t^{2}-2150t+140$&$-6t^{2}+6t$\\\hline
\end{tabular}
\end{adjustbox}
\caption{Trial and test functions with real coefficients over the master element $[0,1]$ that lead to a four-stage Runge-Kutta method.}
\label{solrealRK4}
\end{table}

\pgfmathsetmacro{\width} {5cm}
\pgfmathsetmacro{\height} {5cm}
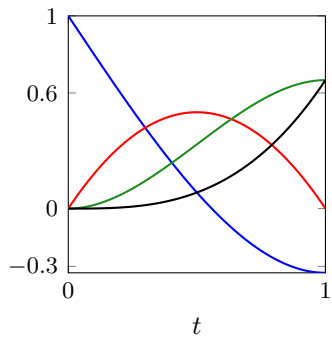
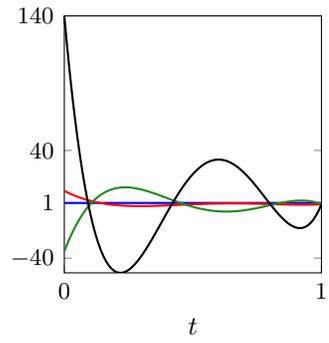
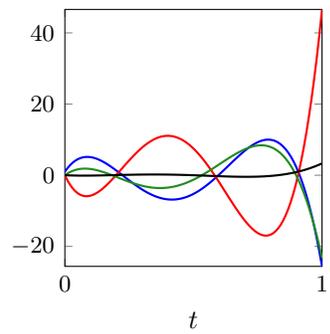
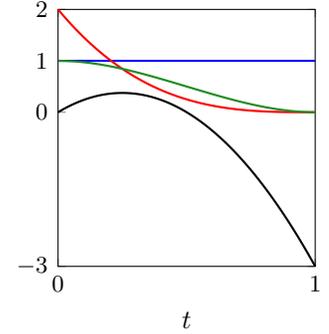
\begin{figure}[h!]
\begin{subfigure}[b]{0.5\linewidth}
\centering
\begin{tikzpicture}
\begin{axis}[footnotesize,
		  width=\width,height=\height, xlabel=$t$,
		 xtick={0,1},ytick={-0.3,0,0.6,1},
		  enlargelimits=false]
\addplot [domain=0:1,samples=100,color=blue,thick]{(2/3)*x*x*x-2*x+1};
\addplot [domain=0:1,samples=100,color=red,thick]{-2*x*x+2*x};
\addplot [domain=0:1,samples=100,color=forestgreen,thick]{-(4/3)*x*x*x+2*x*x};
\addplot [domain=0:1,samples=100,color=black,thick]{(2/3)*x*x*x};
\end{axis}
\end{tikzpicture}
\caption{\footnotesize Trial functions of the first real solution.}  
\end{subfigure}
\begin{subfigure}[b]{0.5\linewidth}
\centering
\begin{tikzpicture}
\begin{axis}[footnotesize,
		  width=\width,height=\height, xlabel=$t$,
		 xtick={0,1},ytick={-40,1,40,140},
		 enlargelimits=false]
\addplot [domain=0:1,samples=100,color=blue,thick]{1};
\addplot [domain=0:1,samples=100,color=red,thick]{140*x*x*x*x-350*x*x*x+300*x*x-100*x+10};
\addplot [domain=0:1,samples=100,color=forestgreen,thick]{-910*x*x*x*x+2170*x*x*x-1725*x*x+500*x-35};
\addplot [domain=0:1,samples=100,color=black,thick]{4480*x*x*x*x-10360*x*x*x+7890*x*x-2150*x+140};
\end{axis}
\end{tikzpicture}
\caption{\footnotesize Test functions of the first real solution.} 
\end{subfigure}
\begin{subfigure}[b]{0.5\linewidth}
\centering
\begin{tikzpicture}
\begin{axis}[footnotesize,
		  width=\width,height=\height, xlabel=$t$,
		  xtick={0,1},ytick={-20,0,20,40},
		  enlargelimits=false]
\addplot [domain=0:1,samples=100,color=blue,thick]{-(2800/3)*x*x*x*x+1610*x*x*x-810*x*x+(320/3)*x+1};
\addplot [domain=0:1,samples=100,color=red,thick]{(4480/3)*x*x*x*x-2520*x*x*x+1230*x*x-(470/3)*x};
\addplot [domain=0:1,samples=100,color=forestgreen,thick]{-(1820/3)*x*x*x*x+980*x*x*x-450*x*x+(160/3)*x};
\addplot [domain=0:1,samples=100,color=black,thick]{(140/3)*x*x*x*x-70*x*x*x+30*x*x-(10/3)*x};
\end{axis}
\end{tikzpicture}
\caption{\footnotesize Trial functions of the second real solution.} 
\end{subfigure}
\begin{subfigure}[b]{0.5\linewidth}
\centering
\begin{tikzpicture}
\begin{axis}[footnotesize,
		  width=\width,height=\height, xlabel=$t$,
		  xtick={0,1},ytick={-3,0,1,2},
		  enlargelimits=false]
\addplot [domain=0:1,samples=100,color=blue,thick]{1};
\addplot [domain=0:1,samples=100,color=red,thick]{-2*x*x*x+6*x*x-6*x+2};
\addplot [domain=0:1,samples=100,color=forestgreen,thick]{2*x*x*x-3*x*x+1};
\addplot [domain=0:1,samples=100,color=black,thick]{-6*x*x+3*x};
\end{axis}
\end{tikzpicture}
\caption{\footnotesize Test functions of the second real solution.} 
\end{subfigure}
\caption{Trial and test functions over the master element $[0,1]$ for the real solutions of the four-stage Runge-Kutta method.}
\label{RK4}
\end{figure}

\clearpage
\section{Conclusions}\label{S4}
We propose a discontinuous Petrov-Galerkin formulation of the linear diffusion equation that lead to explicit Runge Kutta methods. We define families of piecewise polynomials for trial and test functions for any stage Runge-Kutta method. We provide explicit examples for Runge-Kutta methods of up to four stages. When the trial functions are polynomials of order $p$ in time, then the test space is formed by functions of order $p+1$ but incomplete. Alternatively, we can define the test space to be a complete polynomial space of order p and the resulting trial space  will be of incomplete polynomials of order $p+1$. Methods with more than two stages result in systems which have solutions with complex coefficients for the polynomial basis which lead to equivalent Runge-Kutta methods of the appropriate order. A limitation of our construction of variational formulations equivalent to explicit Runge-Kutta methods is that, for a large number of stages, we end up with large nonlinear systems of equations that are difficult to solve.

Possible extensions of this work are to study the variational structure of other implicit and explicit methods such as Adams-Bashforth, Adams-Moulton or Backward Differentiation Formulas (BDF). Additionally, we will apply these formulations to design explicit (goal-oriented) adaptive strategies. Future work, will also analyze the stability of the new time marching schemes arising from our Galerkin construction in order to build more stable explicit methods. 

\appendix

\section{Matrix form of the nonlinear system}\label{A1}
In this section we express (\ref{CondfinalRK}) in matrix form. We consider, for example, $s$-stages and trial and test functions of order $s$ over the master element $[0,1]$ 
\begin{equation}\nonumber
\displaystyle{\left\{
\begin{split}
\phi_{1}(t)&=c_{10}+c_{11}t+\ldots+c_{1s}t^{s},\\
&\;\;\vdots\\
\phi_{s}(t)&=c_{s0}+c_{s1}t+\ldots+c_{ss}t^{s},\\
\varphi_{1}(t)&=d_{10}+d_{11}t+\ldots+d_{1s}t^{s},\\
&\;\;\vdots\\
\varphi_{s}(t)&=d_{s0}+d_{s1}t+\ldots+d_{ss}t^{s}.\\
\end{split}
\right.}
\end{equation}

To simplify notation, we collect the entries into the following matrices
$$C:=\begin{bmatrix}c_{10}&c_{11}&\cdots&c_{1s}\\\vdots&\vdots&\ddots&\vdots\\c_{s0}&c_{s1}&\cdots&c_{ss}\end{bmatrix},$$
$$D:=\begin{bmatrix}d_{10}&d_{11}&\cdots&d_{1s}\\\vdots&\vdots&\ddots&\vdots\\d_{s0}&d_{s1}&\cdots&d_{ss}\end{bmatrix},$$
thus, we can write conditions (\ref{CondfinalRK}) in matrix form as

\begin{equation}\label{sys.app}
\displaystyle{\left\{
\begin{split}
&C\,\mathbf{e}_{1,s+1}=\mathbf{e}_{1,s},\\ 
&D\,\mathbf{1}_{s+1}=\mathbf{e}_{1,s},\\ 
&DAC^{T}-BC^{T}=E,\\ 
&DFC^{T}=G,
\end{split}
\right.}
\end{equation}
where

$$\mathbf{e}_{1,s+1}:=\begin{bmatrix}1&0&\cdots&0\end{bmatrix}^{T},\;\;\mathbf{1}_{s+1}:=\begin{bmatrix}1&1&\cdots&1\end{bmatrix}^{T},$$
$$A:=\begin{bmatrix}0&1&\cdots&1&1\\0&1/2&\cdots&(s-1)/s&s/(s+1)\\\vdots&\vdots&\ddots&\vdots&\vdots\\0&1/s&\cdots&1/2&s/(2s-1)\\0&1/(s+1)&\cdots&(s-1)/(2s-1)&1/2\end{bmatrix},$$
$$B:=\begin{bmatrix}1&1&\cdots&1\\0&0&\cdots&0\\\vdots&\vdots&\ddots&\vdots\\0&0&\cdots&0\end{bmatrix},\;\;E:=\begin{bmatrix}-1&0&\cdots&0\\-1&1&\cdots&0\\\vdots&\vdots&\ddots&\vdots\\-1&0&\cdots&1\end{bmatrix},$$
$$F:=\begin{bmatrix}1&1/2&\cdots&1/s&1/(s+1)\\1/2&1/3&\cdots&1/(s+1)&1/(s+2)\\\vdots&\vdots&\ddots&\vdots&\vdots\\1/s&1/(s+1)&\cdots&1/(2s-1)&1/2s\\1/(s+1)&1/(s+2)&\cdots&1/2s&1/(2s+1)\end{bmatrix},$$
$$G:=\begin{bmatrix}b_{1}&b_{2}&\cdots&b_{s}\\a_{21}&a_{22}&\cdots&a_{2s}\\\vdots&\vdots&\ddots&\vdots\\a_{s1}&a_{s2}&\ldots&a_{ss}\end{bmatrix}.$$

We compute the entries in the matrices $A$ and $F$ from
\begin{equation}\nonumber
\begin{split}
A&=\int_{0}^{1}\begin{bmatrix}1\\t\\\vdots\\t^{s-1}\\t^{s}\end{bmatrix}\begin{bmatrix}0&1&\cdots&(s-1)t^{s-2}&st^{s-1}\end{bmatrix}dt=\\
&=\int_{0}^{1}\begin{bmatrix}0&1&\cdots&(s-1)t^{s-2}&st^{s-1}\\0&t&\cdots&(s-1)t^{s-1}&st^{s}\\\vdots&\vdots&\ddots&\vdots&\vdots\\0&t^{s-1}&\cdots&(s-1)t^{2s-3}&st^{2s-2}\\0&t^{s}&\cdots&(s-1)t^{2s-2}&st^{2s-1}\end{bmatrix}dt=\\
&=\begin{bmatrix}0&1&\cdots&1&1\\0&1/2&\cdots&(s-1)/s&s/(s+1)\\\vdots&\vdots&\ddots&\vdots&\vdots\\0&1/s&\cdots&1/2&s/(2s-1)\\0&1/(s+1)&\cdots&(s-1)/(2s-1)&1/2\end{bmatrix},
\end{split}
\end{equation}
\begin{equation}\nonumber
\begin{split}
F&=\int_{0}^{1}\begin{bmatrix}1\\t\\\vdots\\t^{s-1}\\t^{s}\end{bmatrix}\begin{bmatrix}1&t&\cdots&t^{s-1}&t^{s}\end{bmatrix}dt=\\
&=\int_{0}^{1}\begin{bmatrix}1&t&\cdots&t^{s-1}&t^{s}\\t&t^{2}&\cdots&t^{s}&t^{s+1}\\\vdots&\vdots&\ddots&\vdots&\vdots\\t^{s-1}&t^{s}&\cdots&t^{2s-2}&t^{2s-1}\\t^{s}&t^{s+1}&\cdots&t^{2s-1}&t^{2s}\end{bmatrix}dt=\\
&=\begin{bmatrix}1&1/2&\cdots&1/s&1/(s+1)\\1/2&1/3&\cdots&1/(s+1)&1/(s+2)\\\vdots&\vdots&\ddots&\vdots&\vdots\\1/s&1/(s+1)&\cdots&1/(2s-1)&1/2s\\1/(s+1)&1/(s+2)&\cdots&1/2s&1/(2s+1)\end{bmatrix}.
\end{split}
\end{equation}

\lstset{language=Matlab,
    breaklines=true,
    morekeywords={matlab2tikz},
    keywordstyle=\color{blue},
    morekeywords=[2]{1}, keywordstyle=[2]{\color{black}},
    identifierstyle=\color{black},
    commentstyle=\color{mygreen},
    numbers=left,
    numberstyle={\tiny \color{black}},
    numbersep=9pt,
    emph=[1]{alpha},emphstyle=[1]\color{mylilas},   
}
\section{MATLAB Code}\label{A2}
In this section we provide a MATLAB code to solve system (\ref{CondfinalRK}) in matrix form for a general number of stages.

\begin{lstlisting}
%Script to calculate the trial and test functions of Runge-Kutta methods
%Import data
[n_fun,p_trial,p_test,V_triv,V_deriv,V_grad]=data;
%Initialize the solution
S_trial=cell(n_fun,p_trial+1);
S_test=cell(n_fun,p_test+1);
%Write the coefficients 
coef_trial=sym('c%d%d', [n_fun p_trial+1]);
coef_test=sym('d%d%d', [n_fun p_test+1]);
%Write the conditions
[C_triv,C_deriv,C_grad]=conditions(n_fun,coef_trial,coef_test,p_trial,p_test);
%Solve the nonlinear system
[S_trial{1:end},S_test{1:end}]=solve([C_triv C_deriv C_grad]==[V_triv V_deriv V_grad],[coef_trial coef_test]);
\end{lstlisting}
\vspace{0.3cm}
\begin{lstlisting}
function [n_fun,p_trial,p_test,V_triv,V_deriv,V_grad] = data()
%Function to write the data
%Number of trial and test functions and their order
n_fun=2;
p_trial=2;
p_test=2;
%The value of the conditions. For example: 
%2-stage Runge-Kutta method
V_triv=[1 1 ;0 0];
V_deriv=[-1 0;-1 1];
alpha=1;
V_grad=[1-1/(2*alpha) 1/(2*alpha);alpha 0];
%Forward Euler method
%V_triv=[1 1];
%V_deriv=-1;
%V_grad=1;
end

\end{lstlisting}
\vspace{0.3cm}
\begin{lstlisting}
function [C_triv,C_deriv,C_grad]= conditions(n_fun,coef_trial,coef_test,p_trial,p_test)
%Function to write the conditions in matrix form the trial and test functions must satisfy
C_trial=coef_trial*eye(p_trial+1,1);
C_test=coef_test*ones(p_test+1,1);
C_triv=[C_trial C_test];

C=repmat((1:p_test+1)',1,p_trial+1)+repmat(0:p_trial,p_test+1,1);
C=1./C;
C_grad=coef_test*C*coef_trial.';

B=[ones(1,p_trial+1);zeros(n_fun-1,p_trial+1)];
A=repmat(0:p_trial,p_test+1,1).*[zeros(p_test+1,1) C(:,1:end-1)];
C_deriv=coef_test*A*coef_trial.'-B*coef_trial.';
end

\end{lstlisting}
\vspace{0.3cm}

\section*{Acknowledgments}
All authors have received funding from the European Union's Horizon 2020 research and innovation programme under the Marie Sklodowska-Curie grant agreement No 644602.

David Pardo, Elisabete Alberdi and Judit Mu\~noz-Matute were partially funded by the Basque Government Consolidated Research Group Grant IT649-13 on ``Mathematical Modeling, Simulation, and Industrial Applications (M2SI)'' and the Projects of the Spanish Ministry of Economy and Competitiveness with reference MTM2016-76329-R (AEI/FEDER, EU), and MTM2016-81697-ERC.

David Pardo has also received funding from the BCAM ``Severo Ochoa'' accreditation of excellence SEV-2013-0323 and the Basque Government through the BERC 2014-2017 program.

Victor M. Calo was partially funded by the CSIRO Professorial Chair in Computational Geoscience at Curtin University, the Mega-grant of the Russian Federation Government (N 14.Y26.31.0013) and the Deep Earth Imaging Enterprise Future Science Platforms of the Commonwealth Scientific Industrial Research Organisation, CSIRO, of Australia. Additional support was provided at Curtin University by The Institute for Geoscience Research (TIGeR) and by the Curtin Institute for Computation.

Judit Mu\~noz-Matute has received funding from the University of the Basque Country (UPV/EHU) grant No. PIF15/346.

\section*{References}
\bibliographystyle{elsarticle-num} 
\bibliography{references.bib}

\end{document}